\documentclass[onefignum,onetabnum]{siamart190516}

\usepackage[utf8]{inputenc}
\usepackage{enumitem}
\usepackage{amsmath,amssymb,amsfonts,xcolor}
\usepackage{thmtools, thm-restate}
\usepackage{algorithmic, algorithm}
\crefname{subsection}{subsection}{subsections}
\usepackage{epstopdf}
\ifpdf
  \DeclareGraphicsExtensions{.eps,.pdf,.png,.jpg}
\else
  \DeclareGraphicsExtensions{.eps}
\fi

\usepackage{graphicx,subcaption}
\usepackage{pgfplots}
\pgfplotsset{compat=1.12}
\usepackage{tikz}
\usetikzlibrary{shapes.geometric, arrows}
\usetikzlibrary{shapes.misc}
\usetikzlibrary{positioning}
\usetikzlibrary{calc}
\usetikzlibrary{matrix}
\usetikzlibrary{patterns}
\usetikzlibrary{trees}
\tikzstyle{seps} = [circle, minimum size = 1cm, draw = black, fill=none]
\tikzstyle{arrow} = [thick,->,>=stealth]
\tikzstyle{darrow} = [thick,<->,>=stealth]
\newdimen\nodeDist
\tikzset{cross/.style={cross out, draw=black, fill=none, minimum size=2*(#1-\pgflinewidth), inner sep=0pt, outer sep=0pt}, cross/.default={3pt}}

\usepackage{booktabs}
\usepackage{multirow}
\usepackage{tabularx}

\newcommand{\red}[1]{\textcolor{red}{#1}}

\headers{Hierarchical Orthogonal Factorization: Sparse Least Squares}{A. Gnanasekaran, E. Darve}

\title{Hierarchical Orthogonal Factorization: Sparse Least Squares problems}

\author{Abeynaya Gnanasekaran\thanks{Institute for Computational and Mathematical Engineering, Stanford University, CA
  (\email{abeynaya@stanford.edu}), (\email{darve@stanford.edu})}
\and Eric Darve\footnotemark[1]}

\usepackage{amsopn}

\usepackage{hyperref}
\usepackage{bbm}

\usepackage{microtype}

\usepackage{mathtools}

\begin{document}

\maketitle

\begin{abstract}
    In this work, we develop a fast hierarchical solver for solving large, sparse least squares problems. We build upon the algorithm, spaQR (sparsified QR~\cite{gnanasekaran2020hierarchical}), that was developed by the authors to solve large sparse linear systems. Our algorithm is built on top of a Nested Dissection based multifrontal QR approach. We use low-rank approximations on the frontal matrices to sparsify the vertex separators at every level in the elimination tree. Using a two-step sparsification scheme, we reduce the number of columns and maintain the ratio of rows to columns in each front without introducing any additional fill-in. With this improvised scheme, we show that the runtime of the algorithm scales as $\mathcal{O}(M \log N)$ and uses $\mathcal{O}(M)$ memory to store the factorization. This is achieved at the expense of a small and controllable approximation error. The end result is an approximate factorization of the matrix stored as a sequence of sparse orthogonal and upper-triangular factors and hence easy to apply/solve with a vector. Finally, we compare the performance of the spaQR algorithm in solving sparse least squares problems with direct multifrontal QR and CGLS iterative method with a standard diagonal preconditioner. 
\end{abstract}

\section{Introduction}

In this work, we are interested in solving large sparse linear least squares problems of the form, 
\begin{align}\label{eqn: ls}
\min_x \|Ax-b\|_2, \quad A \in \mathbb{R}^{M \times N},\; M \geq N,\; \text{rank}(A)=N 
\end{align}
Least squares problems appear in a wide range of topics ranging from constrained optimization, computer graphics to economics, data analysis and machine learning. These typically have millions of rows and columns but in many applications may have less than $0.1\%$ of the entries that are non-zero. It is therefore important to maintain the sparsity of the matrix for efficient solutions to large sparse least squares problems. 

The most commonly used technique to solve the least squares problems is to solve the normal equations,  
\begin{align}\label{eqn: normal_eqn}
    A^T A x = A^T b
\end{align}
The normal equations form a symmetric positive definite linear system which can be solved with a sparse Cholesky factorization or an iterative method like CG~\cite{Hestenes&Stiefel:1952}. However, the condition number of the coefficient matrix $A^TA$ in the normal equations is square of the condition number of $A$. Hence, explicitly forming $A^TA$ and solving the normal equations may lead to numerical instabilities. As a result, this method is not appropriate for ill-conditioned or stiff problems~\cite{Bjorck:1411947,matrix_computations}.  

A more reliable and accurate direct method for least squares problems is using the QR decomposition of $A$. Then, solving \Cref{eqn: ls} boils down to the solution of the linear system,
\[
    Rx = Q^Tb
\]
The $R$ factor in the QR decomposition is same as the Cholesky factor of $A^TA$. However, the use of orthogonal transformations in the QR decomposition make the factorization stable as opposed to computing the Cholesky decomposition of $A^TA$. 

Direct methods can be expensive even for sparse matrices due to the fill-in (new non-zero entries) introduced during the factorization. On the other hand, iterative methods such as CGLS~\cite{Golub1965NumericalMF} (also known as CGNR), LSQR~\cite{journals/toms/PaigeS82}, LSMR~\cite{journals/siamsc/FongS11}, which have been proposed for least-squares problems, rarely work well without good preconditioners. A good hybrid between the two are incomplete factorizations, which are then used as preconditioners with iterative methods. For example, Incomplete Cholesky (on $A^TA$)~\cite{Manteuffel1980AnIF}, Incomplete QR (IQR)~\cite{jennings,Saad1988PreconditioningTF}, Multilevel Incomplete QR (MIQR)~\cite{miqr} find an approximate factorization by thresholding the fill-in entries. In general, these preconditioners are not guaranteed to work and can fail for a large number of problems~\cite{James1990ConjugateGM,Chow1997ExperimentalSO,iqr_fail}. Better preconditioners can be built when additional information on the problem is available. 

Hierarchical solvers such as Hierarchical Interpolative Factorization (HIF)~\cite{ho2013hierarchical,Ho2016HierarchicalIF,FeliuFab2018RecursivelyPH,feliufaba2020hierarchical}, LoRaSp~\cite{lorasp1,lorasp2} and Sparsified Nested Dissection (spaND)~\cite{2019arXiv190102971C,klockiewicz2020second} are another family of incomplete factorizations. All three solvers were developed to perform a fast Cholesky factorization of SPD matrices by incorporating low-rank approximations in the classical multifrontal approach. Some of these algorithms were later extended to perform a fast LU factorization on unsymmetric matrices~\cite{ho2013hierarchical}. In contrast, another popular approach is to store the dense fronts using low-rank bases and perform fast matrix-vector products using these bases~\cite{mumps,H_QR,C:LaBRI::CIMI15,Ghysels2016AnEM,blr_pastix,Schmitz2012AFD,Xia2013EfficientSM,Xia2009SuperfastMM}. Most of the efforts in this area have been focused on solving linear systems and has recently been extended to solve Toeplitz least squares problems~\cite{Xi2014SuperfastAS}.

The spaQR algorithm developed by the authors in ~\cite{gnanasekaran2020hierarchical} is another hierarchical solver that builds on the idea that certain off-diagonal blocks in $A$ and $A^TA$ are low-rank. This property was leveraged to build a fast approximate sparse QR factorization of $A$ to solve linear systems in $\mathcal{O}(N \log N)$ time~\cite{gnanasekaran2020hierarchical}. 

\textbf{However, if the spaQR algorithm is used as such without any modifications for tall, thin matrices (least squares problems), we will not have a near linear time algorithm}. For example, we can end up with a final block that has few columns but $O(M-N)$ rows. This would lead to a high computational cost. \textbf{So controlling the aspect ratio (i.e., the ratio number of rows / number of columns) of the diagonal blocks is a key step in our new algorithm.} This is one of the goals of the algorithm in this manuscript. 

The (new) spaQR algorithm produces a sparse orthogonal factorization of $A \in \mathbb{R}^{M \times N}$ in $\mathcal{O}(M \log N)$ time, such that, 
\[
A \approx QW = \prod_i Q_i \prod_j W_j 
\]
where each $Q_i$ is a sparse orthogonal matrix and $W_j$ is either a sparse orthogonal or sparse upper triangular matrix. The orthogonal factor $Q$ is not stored to save on the available memory. The least squares solution is computed by using the corrected seminormal equations (CSNE) approach of~\cite{bjorck_article}. This involves solving,
\[
    W^TW {x} = A^Tb 
\]
along with iterative refinement. $W$ is stored as a sequence of sparse Householder vectors and sparse triangular factors which makes it fast to solve with a vector.



Finally, note that while the $W$ factor is not strictly upper triangular, we still use the term ``fast QR solver" as the algorithm is built on top of the multifrontal QR method. To the best of our knowledge, this is the first work of its kind to combine the ideas of low-rank approximations and QR factorization to build a fast sparse least squares solver.


\subsection{Contribution} 

We extend the spaQR algorithm developed in~\cite{gnanasekaran2020hierarchical} to solve linear least squares problem. Our main contributions are as follows:
\begin{itemize}
    \item We introduce a way to reorder and sparsify the rows in addition to the columns in the sparse multifrontal QR factorization for tall, thin matrices.
    \item Our algorithm keeps the aspect ratio of the diagonal blocks bounded leading to a $\mathcal{O}(M \log N)$ algorithm. 
    \item We show numerical benchmarks on least squares problems arising in PDE constrained optimization problems and on matrices taken from the Suite Sparse Matrix Collection~\cite{suitesparse}.
    \item The C++ code for the implementation is freely available for download and use\footnote{\texttt{https://github.com/Abeynaya/spaQR\_public}}. The benchmarks can be reproduced by using the scripts in the repository. 
\end{itemize}

\paragraph{Organization of the paper} In \Cref{Sec: Algo}, we review the spaQR algorithm for solving linear systems and extend it solve least squares problems. We discuss the complexity of the algorithm in \Cref{sec:complexity}. Finally, in ~\Cref{sec: benchmarks}, we provide numerical results on benchmark problems. 

\section{Algorithm}
\label{Sec: Algo}

We begin with a brief review of the spaQR algorithm for solving linear systems (square matrices). Then, we explain the improvements needed to handle least squares problems where the matrix is now tall and thin. In particular, we describe a row reordering strategy and additional sparsification steps. 

\subsection{Sparsified QR (spaQR) for square matrices}
\label{sec: spaQR_intro}

The spaQR algorithm was proposed by the authors in~\cite{gnanasekaran2020hierarchical} to perform an approximate QR factorization of a sparse square matrix in $\mathcal{O}(N \log N)$ time. The algorithm is built on top of a sparse multifrontal QR factorization. The key step is to continually decrease the size of the vertex separators (obtained using Nested Dissection) by using a low-rank approximation of its off-diagonal blocks. Here, we go over the algorithm and introduce key terminology that will be used in the rest of the paper.

Multifrontal methods are organized as a sequence of factorizations of small dense matrices. The order of factorization is given by the elimination tree and the factorization proceeds from the leaves to the root of the tree. Prior to the factorization, the matrix is typically reordered to minimize the fill-in (new non-zero entries) that appear during the factorization. We use a Nested Dissection ordering for this purpose as it typically leads to minimal fill-in~\cite{George1973NestedDO}.

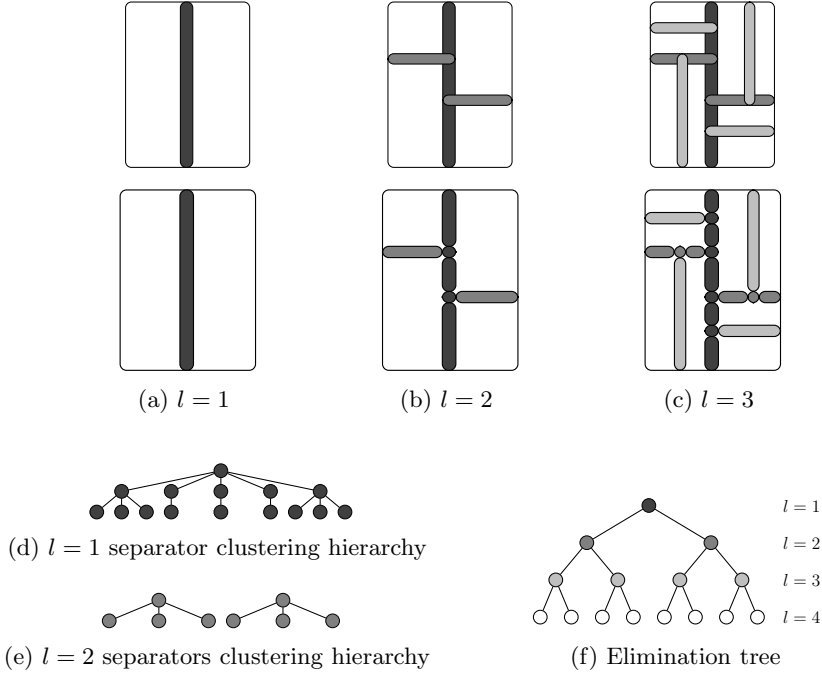
\begin{figure}[tbhp]
    \centering
    \begin{subfigure}[t]{0.25\textwidth}
        \centering
        \scalebox{0.55}{
        \begin{tikzpicture}
            \filldraw[fill = white, rounded corners, line width = 0.25mm] (0, 0) rectangle (3, 4) {};
            \filldraw[fill = darkgray, rounded corners, line width= 0.25mm] (1.325, 0) rectangle (1.625, 4) {};
        \end{tikzpicture}}
    \end{subfigure}%
   ~
    \begin{subfigure}[t]{0.25\textwidth}
        \centering
        \scalebox{0.55}{
        \begin{tikzpicture}
            \filldraw[fill = white, rounded corners, line width = 0.25mm] (0, 0) rectangle (3, 4) {};
            \filldraw[fill = darkgray, rounded corners, line width= 0.25mm] (1.325, 0) rectangle (1.625, 4) {};
            \filldraw[fill= gray, line width =0.25mm, rounded corners] (0, 2.5) rectangle (1.625, 2.75) {};
            \filldraw[fill= gray, line width =0.25mm, rounded corners] (1.325, 1.5) rectangle (3, 1.75) {};
        \end{tikzpicture}}
    \end{subfigure}%
   ~
    \begin{subfigure}[t]{0.25\textwidth}
        \centering
        \scalebox{0.55}{
        \begin{tikzpicture}
            \filldraw[fill = white, rounded corners, line width = 0.25mm] (0, 0) rectangle (3, 4) {};
            \filldraw[fill = darkgray, rounded corners, line width= 0.25mm] (1.325, 0) rectangle (1.625, 4) {};
            \filldraw[fill= gray, line width =0.25mm, rounded corners] (0, 2.5) rectangle (1.625, 2.75) {};
            \filldraw[fill= gray, line width =0.25mm, rounded corners] (1.325, 1.5) rectangle (3, 1.75) {};
            \filldraw[fill= lightgray, line width =0.25mm, rounded corners] (0.65, 0) rectangle (0.9, 2.75) {};
            \filldraw[fill= lightgray, line width =0.25mm, rounded corners] (0, 3.25) rectangle (1.625, 3.5) {};
            \filldraw[fill= lightgray, line width =0.25mm, rounded corners] (2.275, 1.5) rectangle (2.525, 4) {};
            \filldraw[fill= lightgray, line width =0.25mm, rounded corners] (1.325, 0.75) rectangle (3, 1) {};
        \end{tikzpicture}}
    \end{subfigure}%
    \vspace{0.25cm}
    
    \begin{subfigure}[t]{0.25\textwidth}
        \centering
        \scalebox{0.6}{
        \begin{tikzpicture}
            \filldraw[fill = white, rounded corners, line width = 0.25mm] (0, 0) rectangle (3, 4) {};
            \filldraw[fill = darkgray, rounded corners, line width= 0.25mm] (1.325, 0) rectangle (1.625, 4) {};
        \end{tikzpicture}}
        \caption{$l=1$}
    \end{subfigure}%
   ~
    \begin{subfigure}[t]{0.25\textwidth}
        \centering
        \scalebox{0.6}{
        \begin{tikzpicture}
            \filldraw[fill = white, rounded corners, line width = 0.25mm] (0, 0) rectangle (3, 4) {};
            \filldraw[fill = darkgray, rounded corners, line width= 0.25mm] (1.325, 0) rectangle (1.625, 1.5) {};
            \filldraw[fill = darkgray, rounded corners, line width= 0.25mm] (1.325, 1.5) rectangle (1.625, 1.75) {};
            \filldraw[fill = darkgray, rounded corners, line width= 0.25mm] (1.325, 1.75) rectangle (1.625, 2.5) {};
            \filldraw[fill = darkgray, rounded corners, line width= 0.25mm] (1.325, 2.5) rectangle (1.625, 2.75) {};
            \filldraw[fill = darkgray, rounded corners, line width= 0.25mm] (1.325, 2.75) rectangle (1.625, 4) {};
            \filldraw[fill= gray, line width =0.25mm, rounded corners] (0, 2.5) rectangle (1.325, 2.75) {};
            \filldraw[fill= gray, line width =0.25mm, rounded corners] (1.625, 1.5) rectangle (3, 1.75) {};
        \end{tikzpicture}}
        \caption{$l=2$}
        \label{fig: gen_interfaces_lvl2}
    \end{subfigure}%
   ~
    \begin{subfigure}[t]{0.25\textwidth}
        \centering
        \scalebox{0.6}{
        \begin{tikzpicture}
            \filldraw[fill = white, rounded corners, line width = 0.25mm] (0, 0) rectangle (3, 4) {};
           \filldraw[fill = darkgray, rounded corners, line width= 0.25mm] (1.325, 0) rectangle (1.625, 0.75) {};
           \filldraw[fill = darkgray, rounded corners, line width= 0.25mm] (1.325, 0.75) rectangle (1.625, 1) {};
           \filldraw[fill = darkgray, rounded corners, line width= 0.25mm] (1.325, 1) rectangle (1.625, 1.5) {};
            \filldraw[fill = darkgray, rounded corners, line width= 0.25mm] (1.325, 1.5) rectangle (1.625, 1.75) {};
            \filldraw[fill = darkgray, rounded corners, line width= 0.25mm] (1.325, 1.75) rectangle (1.625, 2.5) {};
            \filldraw[fill = darkgray, rounded corners, line width= 0.25mm] (1.325, 2.5) rectangle (1.625, 2.75) {};
            \filldraw[fill = darkgray, rounded corners, line width= 0.25mm] (1.325, 2.75) rectangle (1.625, 3.25) {};
            \filldraw[fill = darkgray, rounded corners, line width= 0.25mm] (1.325, 3.25) rectangle (1.625, 3.5) {};
            \filldraw[fill = darkgray, rounded corners, line width= 0.25mm] (1.325, 3.5) rectangle (1.625, 4) {};
            \filldraw[fill= gray, line width =0.25mm, rounded corners] (0, 2.5) rectangle (0.65, 2.75) {};
            \filldraw[fill= gray, line width =0.25mm, rounded corners] (0.65, 2.5) rectangle (0.9, 2.75) {};
            \filldraw[fill= gray, line width =0.25mm, rounded corners] (0.9, 2.5) rectangle (1.325, 2.75) {};
            \filldraw[fill= gray, line width =0.25mm, rounded corners] (1.625, 1.5) rectangle (2.275, 1.75) {};
            \filldraw[fill= gray, line width =0.25mm, rounded corners] (2.275, 1.5) rectangle (2.525, 1.75) {};
            \filldraw[fill= gray, line width =0.25mm, rounded corners] (2.525, 1.5) rectangle (3, 1.75) {};
            \filldraw[fill= lightgray, line width =0.25mm, rounded corners] (0.65, 0) rectangle (0.9, 2.5) {};
            \filldraw[fill= lightgray, line width =0.25mm, rounded corners] (0, 3.25) rectangle (1.325, 3.5) {};
            \filldraw[fill= lightgray, line width =0.25mm, rounded corners] (2.275, 1.75) rectangle (2.525, 4) {};
            \filldraw[fill= lightgray, line width =0.25mm, rounded corners] (1.625, 0.75) rectangle (3, 1) {};
        \end{tikzpicture}}
        \caption{$l=3$}
        \label{fig: gen_interfaces}
    \end{subfigure}
    
    \vspace{0.25cm}
    
    \begin{subfigure}[b]{0.45\textwidth}
    \centering
    \begin{subfigure}[b]{1\textwidth}
    \centering
    \scalebox{0.55}{
        \begin{tikzpicture}[level distance=0.5cm,
  level 1/.style={sibling distance=1.2cm},
  level 2/.style={sibling distance=0.6cm},
  edge from parent/.append style = {line width = 0.25mm}]
  \node (Root) [circle, draw=black, fill=darkgray] at (0,0) {}
    child {node [circle, draw=black, fill=darkgray]{}
        child {node [circle, draw=black, fill=darkgray]{}}
        child {node [circle, draw=black, fill=darkgray]{}}
        child {node [circle, draw=black, fill=darkgray]{}}
    }
    child {node [circle, draw=black, fill=darkgray]{}
        child {node [circle, draw=black, fill=darkgray]{}}
        }
    child {node [circle, draw=black, fill=darkgray]{} 
        child {node [circle, draw=black, fill=darkgray]{}}
        }
    child {node [circle, draw=black, fill=darkgray]{} 
        child {node [circle, draw=black, fill=darkgray]{}}
        }
    child {node [circle, draw=black, fill=darkgray]{}
        child {node [circle, draw=black, fill=darkgray]{}}
        child {node [circle, draw=black, fill=darkgray]{}}
        child {node [circle, draw=black, fill=darkgray]{}}
    };
    \end{tikzpicture}}
    \caption{$l=1$ separator clustering hierarchy}
    \label{lvl1_ch}
    \end{subfigure}
    
    \begin{subfigure}[b]{1\textwidth}
    \centering
    \scalebox{0.55}{
        \begin{tikzpicture}[level distance=0.5cm,
  level 1/.style={sibling distance=1.2cm},
  edge from parent/.append style = {line width = 0.25mm}]
  \node (Root) [circle, draw=black, fill=gray] at (0, -0.5){}
        child {node [circle, draw=black, fill=gray]{}}
        child {node [circle, draw=black, fill=gray]{}}
        child {node [circle, draw=black, fill=gray]{}};
    \node (Root2) [circle, draw=black, fill=gray] at (3, -0.5){}
        child {node [circle, draw=black, fill=gray]{}}
        child {node [circle, draw=black, fill=gray]{}}
        child {node [circle, draw=black, fill=gray]{}};
    \end{tikzpicture}}
    \caption{$l=2$ separators clustering hierarchy}
    \label{lvl2_ch}
    \end{subfigure}
    
    \end{subfigure}%
    ~
    \begin{subfigure}[b]{0.45\textwidth}
        \centering
        \scalebox{0.55}{
        \begin{tikzpicture}[level distance=0.9cm,
  level 1/.style={sibling distance=3cm},
  level 2/.style={sibling distance=1.5cm},
  level 3/.style={sibling distance=0.75cm},
  edge from parent/.append style = {line width = 0.25mm}]
  \node (Root) [circle, draw=black, fill=darkgray] {}
    child {node [circle, draw=black, fill=gray]{}
      child {node [circle, draw=black, fill=lightgray]{}
        child {node [circle, draw=black]{}}
        child {node [circle, draw=black]{}}
      }
      child {node [circle, draw=black, fill=lightgray]{}
        child {node [circle, draw=black]{}}
        child {node [circle, draw=black]{}}
      }
    }
    child {node [circle, draw=black, fill=gray]{}
      child {node [circle, draw=black, fill=lightgray]{}
        child {node [circle, draw=black]{}}
        child {node [circle, draw=black]{}}
      }
      child {node [circle, draw=black, fill=lightgray]{}
        child {node [circle, draw=black]{}}
        child {node [circle, draw=black]{}}
      }
    };
    \begin{scope}[every node/.style={right}]
     \path (Root    -| Root-2-2-2) ++(5mm,0) node {\large$l=1$};
     \path (Root-1  -| Root-2-2-2) ++(5mm,0) node {\large$l=2$};
     \path (Root-1-1-| Root-2-2-2) ++(5mm,0) node {\large$l=3$};
     \path (Root-1-1-1-| Root-2-2-2) ++(5mm,0) node {\large$l=4$};
   \end{scope}

        \end{tikzpicture}}
        \caption{Elimination tree}
        \label{ND_etree}
    \end{subfigure}
    
    \caption{(a), (b), (c) show the creation of separators and interfaces in the modified Nested Dissection algorithm. (d), (e) show the clustering hierarchy of the interfaces within each separator. (f) shows the elimination tree associated with the Nested Dissection partitioning. }
    \label{mnd_multilvl}
\end{figure}

For computational efficiency, we need to block or partition the separators at every level of the tree. This partitioning needs to be done at all levels, not simply at the level at which the separator was computed. For example,  \Cref{fig: gen_interfaces_lvl2,fig: gen_interfaces} show the multilevel partitioning of the top separator.  The top separator is the dark grey separator in the center.

\begin{figure}[htbp]
    \centering
    \begin{subfigure}[t]{0.3\textwidth}
        \centering
        \scalebox{0.6}{
        \begin{tikzpicture}
            \filldraw[fill = white, rounded corners, line width = 0.25mm] (0, 0) rectangle (3, 4) {};
            \filldraw[fill = gray, rounded corners, line width= 0.25mm] (1.325, 0) rectangle (1.625, 4) {};
            \filldraw[fill= lightgray, line width =0.25mm, rounded corners] (0, 1) rectangle (1.325, 1.25) {};
            \filldraw[fill= lightgray, line width=0.25mm, rounded corners] (1.625, 3) rectangle (3, 3.25) {};
        \end{tikzpicture}}
        \caption{vertex separators}
        \label{seps}
    \end{subfigure}
    \begin{subfigure}[t]{0.3\textwidth}
        \centering
        \scalebox{0.6}{
        \begin{tikzpicture}
            \filldraw[fill = white, rounded corners, line width = 0.25mm] (0, 0) rectangle (3, 4) {};
            \filldraw[fill = gray, rounded corners, line width= 0.25mm] (1.325, 0) rectangle (1.625, 1) {};
            \filldraw[fill = gray, rounded corners, line width= 0.25mm] (1.325, 1) rectangle (1.625, 1.25) {};
            \filldraw[fill =gray, rounded corners, line width= 0.25mm] (1.325, 1.25) rectangle (1.625, 3) {};
            \filldraw[fill = gray, rounded corners, line width= 0.25mm] (1.325, 3) rectangle (1.625, 3.25) {};
            \filldraw[fill = gray, rounded corners, line width= 0.25mm] (1.325, 3.25) rectangle (1.625, 4) {};
            \filldraw[fill= lightgray, line width =0.25mm, rounded corners] (0, 1) rectangle (1.325, 1.25) {};
            \filldraw[fill= lightgray, line width=0.25mm, rounded corners] (1.625, 3) rectangle (3, 3.25) {};
        \end{tikzpicture}}
        \caption{interfaces}
        \label{int}
    \end{subfigure}
    \caption{Nested dissection partition with 3 levels. The figure on the left shows the vertex separators and the one on the right shows the interfaces. }
    \label{interfaces_full}
\end{figure}
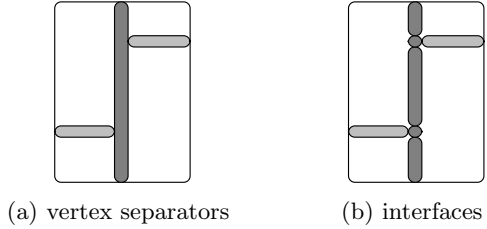

In the rest of the paper, we will call each block in this partition an interface. \Cref{interfaces_full} shows a 3-level Nested Dissection partition of an arbitrary graph; the figure on the left shows the Nested Dissection separators and the one on the right shows the interfaces. To minimize the number of floating point operations and the rank of the off-diagonal blocks, it is best if an interface is connected to few subdomains in the nested dissection. At a minimum, an interface will be connected to two subdomains. But because of the way separators are constructed, we will also have small interfaces that are connected to three or more subdomains. 


The spaQR algorithm alternates between factoring (block QR) the separators at a level $l$, scaling and sparsifying the `interfaces' at all levels $l'>l$ in the elimination tree. Consider the spaQR algorithm with a simple 3-level Nested Dissection partitioning as shown in \Cref{interfaces_full}. The spaQR algorithm first performs a block QR factorization on each interior (white region in \Cref{interfaces_full}). Then, it scales and sparsifies the interfaces belonging to the remaining separators as discussed next. Consider a subset of the top separator in \Cref{interfaces_full} that is at the interface between two interiors. Let $p$ be the interface and $n$ be all the nodes in the subdomains connected to the interface. Then, the associated block in the matrix is, 
\[A_p = 
\begin{bmatrix}
A_{pp} & A_{pn} \\
A_{np} & A_{nn}
\end{bmatrix}
\]
Assume that the block $A_{pp} = I$. This is not a limiting assumption, as we perform a block diagonal scaling on the interfaces (see \Cref{sec: scaling}) before sparsification. The key assumption in the spaQR algorithm is that the off-diagonal blocks ($A_{pn}$, $A_{np}$) are low-rank. Compute a low-rank approximation of,
\[\begin{bmatrix}
A_{np}^T & A_{pn}
\end{bmatrix}  = Q_{pp}W_{pn} = \begin{bmatrix}
Q_{pf} & Q_{pc}
\end{bmatrix}\begin{bmatrix}
W_{fn} \\
W_{cn}
\end{bmatrix} \text{with } \|W_{fn}\|_{_2}=\mathcal{O}(\epsilon)\]
Then apply $Q_{pp}$ to $A$,
\[
\begin{bmatrix}
Q_{pp} &  \\
 & I
\end{bmatrix}^T
\begin{bmatrix}
I & A_{pn} \\
A_{np} & A_{nn}
\end{bmatrix}
\begin{bmatrix}
Q_{pp} &  \\
 & I
\end{bmatrix} = \begin{bmatrix}
I &  & \mathcal{O}(\epsilon) \\
 & I & \hat{A}_{cn} \\
 \mathcal{O}(\epsilon) & \hat{A}_{nc} & A_{nn}
\end{bmatrix}
\]

Thus by applying the orthogonal transformation $Q$, we split the nodes of interface $p$ into `fine' $f$ and `coarse' $c$ nodes. The fine nodes are disconnected from the rest after dropping the $\mathcal{O}(\epsilon)$ terms. Hence, the degree of freedom in interface $p$ has been reduced by $|f|$. The other interfaces can be sparsified similarly. This is the key step in the spaQR algorithm.



\subsection{Partitioning}
\label{ord_clus}

Next, we explain row and column reordering strategy we use to minimize fill-in during the factorization. The column reordering is done in a similar fashion as in ~\cite{gnanasekaran2020hierarchical} for square matrices. However, the row reordering strategy is quite different as least squares problem have more rows than columns. 

The vertex separators are obtained through a Nested Dissection partitioning on the graph of $A^TA$. This can be done either using a graph partitioning software like Metis or a hypergraph based partitioning like PaToH~\cite{atalyrek2011PaToHT}, hMetis~\cite{Karypis1998HmetisAH} and Zoltan~\cite{ZoltanHypergraphIPDPS06}. In our implementation, we provide the option partition using Metis~\cite{metis} and PaToH~\cite{atalyrek2011PaToHT}. Also, for problems where the underlying geometry is available, we can use it to efficiently partition the matrix. The quality of vertex separators and interfaces depend on the partitioning technique. Defining the `best' technique is beyond the scope of this work. 

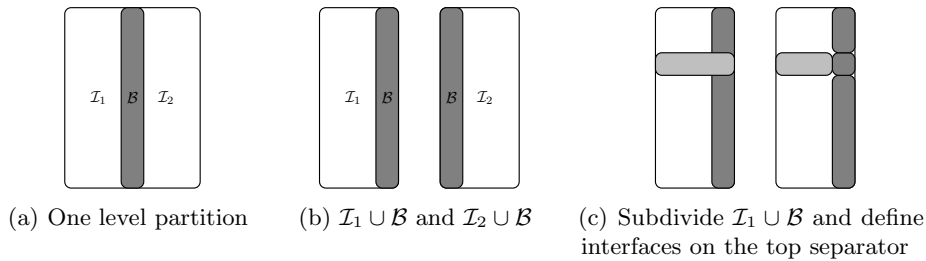
\begin{figure}[!tbhp]
    \centering
    \begin{subfigure}[t]{0.25\textwidth}
        \centering
        \scalebox{0.6}{
        \begin{tikzpicture}
            \filldraw[fill = white, rounded corners, line width = 0.25mm] (0, 0) rectangle (3, 4) {};
            \filldraw[fill = gray, rounded corners, line width= 0.25mm] (1.25, 0) rectangle (1.75, 4) {};
            \node at (0.75,2) {$\mathcal{I}_1$};
            \node at (1.5,2) {$\mathcal{B}$};
            \node at (2.25,2) {$\mathcal{I}_2$};
        \end{tikzpicture}}
        \caption{One level partition}
    \end{subfigure}%
   ~
    \begin{subfigure}[t]{0.3\textwidth}
        \centering
        \scalebox{0.6}{
        \begin{tikzpicture}
            \filldraw[fill = white, rounded corners, line width = 0.25mm] (0, 0) rectangle (1.75, 4) {};
            \filldraw[fill = gray, rounded corners, line width= 0.25mm] (1.25, 0) rectangle (1.75, 4) {};
            \node at (0.75,2) {$\mathcal{I}_1$};
            \node at (1.5,2) {$\mathcal{B}$};
        \end{tikzpicture}
        }
        \quad
        \scalebox{0.6}{
        \begin{tikzpicture}
            \filldraw[fill = white, rounded corners, line width = 0.25mm] (0, 0) rectangle (1.75, 4) {};
            \filldraw[fill = gray, rounded corners, line width= 0.25mm] (0, 0) rectangle (0.5, 4) {};
            \node at (0.25,2) {$\mathcal{B}$};
            \node at (1,2) {$\mathcal{I}_2$};
        \end{tikzpicture}}
        \caption{$\mathcal{I}_1 \cup \mathcal{B}$ and $\mathcal{I}_2 \cup \mathcal{B}$ }
    \end{subfigure}%
   ~
    \begin{subfigure}[t]{0.35\textwidth}
        \centering
        \scalebox{0.6}{
        \begin{tikzpicture}
            \filldraw[fill = white, rounded corners, line width = 0.25mm] (0, 0) rectangle (1.75, 4) {};
            \filldraw[fill = gray, rounded corners, line width= 0.25mm] (1.25, 0) rectangle (1.75, 4) {};
            \filldraw[fill= lightgray, line width =0.25mm, rounded corners] (0, 2.5) rectangle (1.75, 3) {};
        \end{tikzpicture}}
        \quad
        \scalebox{0.6}{
        \begin{tikzpicture}
            \filldraw[fill = white, rounded corners, line width = 0.25mm] (0, 0) rectangle (1.75, 4) {};
            \filldraw[fill = gray, rounded corners, line width= 0.25mm] (1.25, 0) rectangle (1.75, 2.5) {};
            \filldraw[fill = gray, rounded corners, line width= 0.25mm] (1.25, 2.5) rectangle (1.75, 3) {};
            \filldraw[fill =gray, rounded corners, line width= 0.25mm] (1.25, 3) rectangle (1.75, 4) {};
            \filldraw[fill= lightgray, line width =0.25mm, rounded corners] (0, 2.5) rectangle (1.25, 3) {};
        \end{tikzpicture}}
        \caption{Subdivide $\mathcal{I}_1 \cup \mathcal{B}$ and define interfaces on the top separator}
    \end{subfigure}
    \caption{The first figure shows a one level partition of an arbitrary graph (hypergraph) using nested dissection (HUND). The next two figures depict the process of identifying the interfaces by subdividing $\mathcal{I}_1 \cup \mathcal{B}$.}
    \label{mnd}
\end{figure}

Once we obtain the vertex separators, we can cluster the unknowns in the separator to define interfaces. We use the modified version of Nested Dissection described in~\cite{2019arXiv190102971C} to do this. The idea is to keep track of the boundary $\mathcal{B}$ of each interior $\mathcal{I}$ and use this in the dissection process. One step of this scheme is shown in \Cref{mnd}. 
A detailed description of this scheme is given in Algorithm 2.2 of~\cite{2019arXiv190102971C}. The columns of the matrix are reordered following the Nested Dissection ordering and the clustering hierarchy. 

Least squares problems have more rows that columns. We want to assign rows to each column cluster such that the diagonal blocks in the reordered matrix will have a small condition number. Two strategies are used for the same. First, we perform a bipartite matching between the rows and columns of the matrix. This matches every column to a unique row of the matrix. The MC64 routine from the HSL Mathematical Software Library~\cite{hsl_mc64} is used to perform the matching. This leaves us with $M-N$ unassigned rows. Each unassigned row $r$ is assigned to cluster $c$ such that, $c = \arg\max_{c_k} \sum_{j \in c_k} A^2_{rj}$ where the summation is over all the nodes belonging to cluster $c_k$. This heuristic identifies the cluster such that the weight of the row is maximized in that cluster. We found that this heuristic works well in practice.

\subsection{Householder QR on Separators}
\label{sec: houseQR}

The spaQR algorithm alternates between classical interiors/separators factorization and sparsification of interfaces. In this section, we describe the factorization of interiors or separators at a level $l$ using block Householder QR.

Let $s$ be the separator of interest, $n$ be all its neighbors (i.e., $(A^TA)_{ns} \ne 0)$ and $w$ be the rest of the nodes disconnected from $s$ in the graph of $A^TA$. Let nodes in $n$ be further categorized into $n=\{n_1, n_2, n_3 \}$. Nodes $n_1$ are such that $A_{n_1s} \ne 0$, while $A_{sn_1}$ may or may not be zero. Nodes $n_2$ are such that $A_{n_2s} = 0$ and $A_{sn_2} \ne 0$ and nodes $n_3$ are such that $A_{n_1n_3} \ne 0$ and $A_{sn_3}=0$. All such nodes $n$ will correspond to $(A^TA)_{ns} \ne 0$. Let $w$ be the remaining nodes in the matrix. Consider the matrix A blocked in the following form, 
\[ A = \begin{bmatrix}
A_{ss} & A_{sn_1} & A_{sn_2} & & \\
A_{n_1s} & A_{n_1n_1} & & A_{n_1n_3} & \\
 &  &  A_{n_2n_2} & A_{n_2n_3} & A_{n_2w} \\
  & A_{n_3n_1} &  A_{n_3n_2} & A_{n_3n_3} & A_{n_3w} \\
& A_{wn_1} & A_{wn_2} & A_{wn_3} & A_{ww}
\end{bmatrix} \]
When dealing with tall thin matrices, the diagonal blocks corresponding to each separator/interior are also tall and thin. Let us denote the number of rows in a diagonal block $A_{kk}$ as $r_k$ and the number of columns as $c_k$. Then, the block Householder transformation $H$ on the columns of separator $s$ gives us,
\[ H^T  \begin{tikzpicture}[baseline=(current bounding box.center)]
\matrix(M) [matrix of math nodes, nodes in empty cells,
                left delimiter={[}, right delimiter={]},
                row sep={1.2em,between origins},
                minimum width=width("An1s")]
        {\\
        A_{ss}  \\
        \\
        \\
        A_{n_1s}\\
        \\
        };
\draw (M-1-1.north west) rectangle (M-3-1.south east);
\draw (M-4-1.north west) rectangle (M-6-1.south east);
\end{tikzpicture}
 = \begin{tikzpicture}[baseline=(current bounding box.center)]
\matrix(M) [matrix of math nodes, nodes in empty cells,
                left delimiter={[}, right delimiter={]},
                row sep={1.2em,between origins},
                minimum width=width("An1s")]
        { \\
        \hat{R}_{ss}\\
        \\
        \\
        \\
        \\
        };
\draw (M-1-1.north west) rectangle (M-3-1.south east);
\end{tikzpicture} \quad \text{where, } \hat{R}_{ss} = \begin{tikzpicture}[baseline=(current bounding box.center)]
\matrix(M) [matrix of math nodes, nodes in empty cells,
                left delimiter={[}, right delimiter={]},
                column sep={2.4em,between origins},
                row sep={1.2em,between origins}]
        {\\ 
        \\
        \\
        };
\draw[fill] (M-1-1.north west) -- (M-1-1.north east) -- (M-1-1.south east) -- cycle;
\end{tikzpicture}
\] 
Then,
\[
H_s^T A = \begin{bmatrix}
\hat{R}_{ss} & \hat{R}_{sn_1} & \hat{R}_{sn_2} & \hat{R}_{sn_3} & \\
& \hat{A}_{n_1n_1} & \hat{A}_{n_1n_2}& \hat{A}_{n_1n_3} & \\
&  &  A_{n_2n_2} & A_{n_2n_3} & A_{n_2w} \\
  & A_{n_3n_1} &  A_{n_3n_2} & A_{n_3n_3} & A_{n_3w} \\
& A_{wn_1} & A_{wn_2} & A_{wn_3} & A_{ww}
\end{bmatrix} \quad \text{where, } H_s = \begin{bmatrix}
H & &\\
& I_{r_n} & \\
&& I_{r_w}
\end{bmatrix}\]
Note that each $\hat{R}_{sn_i}$ block is rectangular of size $r_s \times c_{n_i}$. Only the top $c_s$ rows go into the final $R$ matrix. In the multifrontal QR, the extra $r_s-c_s$ rows below the main diagonal are added to the parent of $s$ in the elimination tree. In our implementation, we split the extra rows between the clusters comprising $n_1$ (which can include any of the ancestors of $s$ in the elimination tree), following the maximum norm heuristic outlined in \Cref{ord_clus}. Following this permutation, we have,
\[
P_s^T H_s^T A = \begin{bmatrix}
R_{ss} & R_{sn_1} & R_{sn_2} & R_{sn_3} & \\
& \Tilde{A}_{n_1n_1} & \Tilde{A}_{n_1n_2}& \Tilde{A}_{n_1n_3} & \\
&  &  {A}_{n_2n_2} & {A}_{n_2n_3} & {A}_{n_2w} \\
  & {A}_{n_3n_1} &  {A}_{n_3n_2} & {A}_{n_3n_3} & {A}_{n_3w} \\
& A_{wn_1} & A_{wn_2} & A_{wn_3} & A_{ww}
\end{bmatrix} = \begin{bmatrix}
R_{ss} & R_{sn} & \\
& \Tilde{A}_{nn}& {A}_{nw} \\
& A_{wn} & A_{ww}
\end{bmatrix}
\]
where, each $R_{sn_i}$ block corresponds to the first $s$ rows of the block $\hat{R}_{sn_i}$.  The cluster $s$ is now disconnected from the rest. We have only introduced fill-in between the neighbors $n$ and not affected the blocks involving $w$.  Then we proceed by only focusing on the trailing matrix, 

\[P_s^T H_s^T A R_s^{-1} = \begin{bmatrix}
I_s &  & \\
& \Tilde{A}_{nn} & {A}_{nw} \\
& A_{wn} & A_{ww}
\end{bmatrix} \quad \text{ where, } R_s = \begin{bmatrix}
R_{ss} & R_{sn} & \\
& I_{c_n} & \\
& & I_{c_w}
\end{bmatrix} \]

\subsection{Scaling of Interfaces}
\label{sec: scaling}

The spaQR algorithm performs a block diagonal scaling on all the remaining interfaces after every level of separator/interior factorization. Consider an interface $p$, its neighbors $n$, and the associated sub-matrix in $A$,
\[A_p = \begin{bmatrix}
A_{pp} & A_{pn} \\
A_{np} & A_{nn}
\end{bmatrix}\]
The block $A_{pp} \in \mathbb{R}^{r_p \times c_p}$, with $r_p \geq c_p$. Note that, $A_p$ is only a sub-block of the matrix $A$ involving the nodes $p$ and $n$. There can exist other nodes $w$, such that $A_{nw} \neq 0$ or $A_{wn} \neq 0$. However, these blocks are not affected during the scaling and sparsification of $p$ and hence are not considered in the rest of the analysis. 

Find the QR decomposition of $A_{pp}$; $A_{pp} = U_{pp} R_p$, where $R_p = \begin{bmatrix}
 R_{pp} \\
 0
 \end{bmatrix}$, $U_{pp} \in \mathbb{R}^{r_p \times r_p}$, and $R_{pp} \in \mathbb{R}^{c_p \times c_p}$. Then,
 \[U_{p}^TA_pR_p^{-1} = \begin{bmatrix}
 I_{c_p} & \Tilde{A}_{p_1n} \\
 0 & \Tilde{A}_{p_2n}\\
 \Tilde{A}_{np} & A_{nn}
 \end{bmatrix} \quad \text{where, } U_p = \begin{bmatrix}
 U_{pp}  &\\
 & I_{r_n}
 \end{bmatrix} \text{ and } R_p = \begin{bmatrix}
 R_{pp} & \\
 & I_{c_n}
 \end{bmatrix}\]
 The matrix blocks $A_{p_1n}$, $A_{p_2n}$ are such that, $U_{pp}^T A_{np} = \begin{bmatrix}
 \Tilde{A}_{p_1n} \\
 \Tilde{A}_{p_2n} \\
 \end{bmatrix}$ and $\Tilde{A}_{p_1n}$ is $c_p \times c_n$ matrix. Similarly we scale the diagonal blocks belonging to the other interfaces. 
 
The scaling step was shown to improve the accuracy and performance of the preconditioner for solving linear systems  in~\cite{gnanasekaran2020hierarchical}. While the scaling step was optional for linear systems (square matrices), it is necessary for least squares problems. This is because the sparsification step, described next, operates on the two blocks $A_{p_1n}$ and $A_{p_2n}$ individually and not on $A_{pn}$ as a whole.   

\subsection{Sparsification of Interfaces}
\label{sec: sparsify}

Once the separators/interiors at a level $l$ have been factorized and the remaining interfaces scaled, the final step is sparsification. This is different from the interface sparsification described for square matrices in~\cite{gnanasekaran2020hierarchical}. 

The sparsification of the interfaces is done in two step as described next. Consider an interface $p$, its neighbors $n$, and the associated block matrix after block diagonal scaling,
\[A_p = 
\begin{bmatrix}
I_{c_p} & A_{p_1n}   \\
0 & A_{p_2n} \\
A_{np} & A_{nn} 
\end{bmatrix}\]
Assume the off-diagonal blocks $A_{np}$, $A_{p_1n}$ and $A_{p_2n}$ are low-rank.

\subsubsection{Step 1: Compress rows $A_{p_2n}$}
\label{sec: sparsify_step1}

Consider a low rank approximation of $A_{p_2n}$ as follows, 
\[ A_{p_2n} = Q'_{pp}W'_{pn} = \begin{bmatrix}
Q'_{pc} & Q'_{pf}
\end{bmatrix}\begin{bmatrix}
W'_{cn} \\
W'_{fn}
\end{bmatrix} \quad \text{with} \quad \|W'_{fn}\|_{_2}=\mathcal{O}(\epsilon) \]
\[ (Q'_p)^T A_p = \begin{bmatrix}
I_{c_p} & A_{p_1n}   \\
 & W'_{cn}  \\
 & \red{W'_{fn}}  \\
A_{np} & A_{nn} 
\end{bmatrix} \quad \text{where, } Q'_p = \begin{bmatrix}
I_{c_p} &&\\
& Q'_{pp} & \\
& & I_{r_n} 
\end{bmatrix} \]
The block $\|W'_{fn}\|= \mathcal{O}(\epsilon)$ and the corresponding rows can be ignored. This step is not necessary in the spaQR factorization of square matrices. However, for least squares systems, this step keeps the aspect ratio of the diagonal blocks bounded throughout the algorithm. A bounded aspect ratio is important to ensure a linear complexity as we will show in \Cref{sec:complexity}.

\subsubsection{Step 2: Compress block [$A_{np}^T$ \; $A_{p_1n}$]
}
\label{sec: sparsify_step2}

Consider a low rank approximation of, 
\[ \begin{bmatrix}
 A_{np}^T & A_{p_1n}
 \end{bmatrix} = Q_{pp}W_{pn} = \begin{bmatrix}
Q_{pf} & Q_{pc}
\end{bmatrix}\begin{bmatrix}
W_{fn}^{(1)} & W_{fn}^{(2)} \\
W_{cn}^{(1)} & W_{cn}^{(2)}
\end{bmatrix} \quad \text{with, } \|
W_{fn}^{(1)} \; W_{fn}^{(2)}
\|_{_2}=\mathcal{O}(\epsilon) \]
Then, 
\[Q_p^T A_p Q_p = \begin{bmatrix}
 I_f & & \red{W_{fn}^{(2)}}  \\
 & I_c & {W}^{(2)}_{cn} \\
  &  & {W'}_{cn}  \\
 \red{W_{fn}^{(1)T}} & {W}^{(1)T}_{cn} & A_{nn} 
\end{bmatrix} \quad \text{where, } Q_{p} = \begin{bmatrix}
Q_{pp} & \\
& I 
\end{bmatrix} \]
where the Identity block is appropriately sized. We can proceed by dropping the $\mathcal{O}(\epsilon)$ terms and focusing on the trailing matrix,
\[ 
 Q_p^T (Q'_p)^T A_p Q_p = \begin{bmatrix}
 I_f & & \\
 & \Tilde{A}_{cc} & \Tilde{A}_{cn} \\
  & {W}^{(1)T}_{cn} & A_{nn}
\end{bmatrix} \quad \text{where, } \Tilde{A}_{cc} = \begin{bmatrix}
I_{c} \\
0
\end{bmatrix} \text{ and } \Tilde{A}_{cn} = \begin{bmatrix}
W^{(2)}_{cn}\\
W'_{cn}
\end{bmatrix}
\]
In this process, we have disconnected $f$ (`fine' nodes) from the rest and decreased the size of interface $p$ by $|f|$. The $A_{nn}$ and the $A_{nw}$, $A_{ww}$ (not shown here) blocks are not affected. In~\cite{gnanasekaran2020hierarchical}, we give a detailed proof to show that sparsification does not affect the elimination tree of $A^TA$, that is, any two disjoint subtrees remain disjoint after sparsification. This still holds and can be shown by following the exact steps given in Appendix C of~\cite{gnanasekaran2020hierarchical}.

The spaQR for square matrices sparsifies the interfaces by performing only the second step shown here. However, the first step is crucial for tall, thin matrices. In practice, we perform step 1 on all remaining interfaces at a given $l$ and then perform step 2. This is in contrast to performing step 1 and step 2 on a particular interface before moving on to the next interface. This approach is more expensive as the size of $A_{np}$ block for some interface $p$ can be bigger (when the rows below the diagonal of $n$ have not been sparsified yet). 

\subsection{Merging of clusters}
\label{sec:merge}

Once the factorization of separators at a level is done, the interfaces of the remaining ND separators are merged following the cluster hierarchy. For example, in  \Cref{mnd_multilvl}, once the leaves $l=4$ and the $l=3$ separators are factorized, the interfaces of the separators at $l=$ 1, 2 are merged following the clustering hierarchy shown in \Cref{lvl1_ch,lvl2_ch}. Merging simply means combining the block rows and columns of the interfaces into a single block matrix. 

\subsection{Sparsified QR}
\label{sec: algo}

We can now explain the complete modified spaQR algorithm for rectangular (tall, thin) matrices. First, we scale the columns of the matrix, so that the 2-norm of each column is a constant. Next, we obtain vertex separators and interfaces by partitioning the matrix using Metis or PaToH. This gives us a reordering for the columns of the matrix. Finally, we find a good row permutation by following the heuristics outlined in \Cref{ord_clus}. Then, we apply our hierarchical solver, spaQR, to factorize the reordered matrix. The modified spaQR algorithm performs the following sequence of operations at every level $l$ in the elimination tree: block Householder factorizations $(H_s,R_s)$ on the interiors/separators and reassigning the rows below the diagonal of each factorized interior/separator $P_s$ (see \Cref{sec: houseQR}), interface scaling $(U_p, R_p)$ (see \Cref{sec: scaling}), interface sparsification $(Q_p', Q_p)$ (see \Cref{sec: sparsify}) (including some permutations to move the fine nodes and merging of the clusters). This leads to a factorization of the form,
\[ Q^T A W^{-1} \approx I\] 
where,
\begin{align*}
    Q &= \prod_{l=1}^{L}\Bigg( \prod_{s\in S_l}H_s P_s \prod_{p\in C_l}U_p \prod_{p\in C_l} Q_p' \prod_{p\in C_l} Q_p\Bigg)\\
    W &= \prod_{l=L}^{1}\Bigg(\prod_{p\in C_l} Q_p^T\prod_{p\in C_l}R_p \prod_{s\in S_l}R_s  \Bigg)
\end{align*}
\begin{algorithm}
\caption{Sparsified QR (spaQR) algorithm}
\begin{algorithmic}[1]
  \REQUIRE {Sparse matrix A, Tolerance $\epsilon$}
   \STATE {Compute column and row partitioning of A, infer separators and interfaces (see \Cref{ord_clus})}
   \FORALL{$l=L, L-1, \dots 1$}
   \FORALL{separators $s$ at level $l$}
        \STATE {Factorize $s$ using block Householder (see \Cref{sec: houseQR})}
        \STATE {Append $H_s$ to $Q$ and $R_s$ to $W$}
        \STATE {Reassign the rows in $s$ below the diagonal to other neighbor interfaces (see \Cref{sec: houseQR})}
        \STATE {Append $P_s$ to $Q$ }
        
   \ENDFOR
   \FORALL{interfaces $p$ remaining at level $l$}
        \STATE {Perform block diagonal scaling on  $p$ (see \Cref{sec: scaling})}
        \STATE{Append $U_p$ to $Q$ and $R_p$ to $W$}
   \ENDFOR
   \FORALL{interfaces $p$ remaining at level $l$}
        \STATE {Sparsify interface $p$ by performing step 1 of the sparsification process (see \Cref{sec: sparsify_step1})}
        \STATE{Append $Q_p'$ to $Q$}
   \ENDFOR
   \FORALL{interfaces $p$ remaining at level $l$}
        \STATE {Sparsify interface $p$ by performing step 2 of the sparsification process  (see \Cref{sec: sparsify_step2})}
        \STATE{Append $Q_p$ to $Q$ and $Q_p^T$ to $W$}
   \ENDFOR
   \FORALL{separators $s$ remaining at level $l$}
   \STATE {Merge interfaces of $s$ one level following the cluster hierarchy (see \Cref{sec:merge})}
   \ENDFOR
   \ENDFOR
  \RETURN {$Q = \prod_{l=1}^{L}\Bigg( \prod_{s\in S_l}H_s P_s \prod_{p\in C_l}U_p \prod_{p\in C_l} Q_p'\prod_{p\in C_l} Q_p\Bigg)$ \\
  \qquad \qquad $W = \prod_{l=L}^{1}\Bigg(\prod_{p\in C_l} Q_p^T\prod_{p\in C_l}R_p \prod_{s\in S_l}R_s  \Bigg)$ 
 such that $Q^TAW^{-1} \approx I$ }
\end{algorithmic}
\label{Algo: spaQR}
\end{algorithm}

Here, $S_l$ is the set of all separators at level $l$ in the elimination tree and $C_l$ is the set of all interfaces remaining after factorization of separators at level $l$. $Q$ is a product of orthogonal matrices and $W$ is a product of upper triangular and orthogonal matrices. Since, $Q$ and $W$ are available as sequence of elementary transformations, they are easy to invert. The complete algorithm is presented in \Cref{Algo: spaQR}.

\section{Complexity analysis}
\label{sec:complexity}

Let us now discuss the complexity of spaQR for tall, thin matrices. Consider the Nested Dissection process on the graph of $A^TA$ ($G_{A^TA}$). Let us define a node as a subgraph of $G_{A^TA}$.The root of the tree corresponds to $l=1$ and the root node is the entire graph $G_{A^TA}$. The children nodes are subgraphs of $G_{A^TA}$ disconnected by a separator. Let us assume the following properties on the matrices,
\begin{enumerate}
    \item The leaf nodes in the elimination tree contain at most $N_0 \in \mathcal{O}(1)$ nodes.
    \item Let $D_i$ be the set of all nodes $j$ that are descendants of a node $i$, whose size is at least $n_i/2$. We assume that the size of $D_i$ is bounded, i.e, $|D_i|=\mathcal{O}(1) \; \forall i$.
    \item All the Nested Dissection separators are minimal. That is, every vertex in the separator is connected to two disconnected nodes in $G_{A^TA}$.
    \item The number of edges leaving a node (subgraph) of size $n_i$ is at most $n_i^{2/3}$. In other words, a node of size $n_i$ is connected to at most $n_i^{2/3}$  vertices in $G_{A^TA}$. Most matrices that arise in the discretization of 2D and 3D PDEs satisfy this property.
\end{enumerate}  

\paragraph{Multifrontal QR} We first estimate the cost of multifrontal QR on sparse matrices under our assumptions. Consider matrix A of size $M \times N$, $\alpha = M/N$, obtained from a PDE discretization on a 3D grid (for example, matrices obtained from PDE constained optimization problems as shown in \Cref{subsec: 3d_invpoi}). The node at the top of the elimination tree has a size of $N$. By assumption 4, the associated separator has size of, 
\[
c_{\text{top}} \in \mathcal{O}(N^{2/3})
\]
The cost of Householder QR on the block corresponding to the last separator is $\mathcal{O}(r_l c_l^2)$ where $r_l$ is the number of rows in the separator block.  At the time of factorization of the top separator, 
\[
r_{\text{top}} \in \mathcal{O}\big((\alpha-1)N + N^{2/3}\big)
\]
as all the extra rows of other separators get moved to the end. Then, the cost of factorizing the top separator block is, 
\[
h_{\text{top}} = \mathcal{O}\big((\alpha -1)N^{7/3}+N^2\big)
\]
The total cost of the multifrontal QR on the matrix is at least $h_{\text{top}}$
\[
    t_{\text{QR, fact}} = \Omega\big((\alpha -1)N^{7/3}+N^2\big) \quad
    \text{(the notation $\Omega$ means bounded from below)}
\]

Next, we estimate the cost of applying the factorization to a vector $b \in \mathbb{R}^{M}$. We are using the CSNE aproach~\cite{bjorck_article}, which involves solving $R^T R x = A^T b$. Thus, we need to evaluate the cost of solving a linear system with $R$ and $R^T$. A node $i$ at a level $l$ in the elimination tree is of size $2^{-l+1}N \leq n_i \leq 2^{-l+2}N$. The associated separator has a size of at most $\mathcal{O}(n_i^{2/3})$ and $\mathcal{O}(n_i^{2/3})$ non-zeros per row. Then the cost of solving a linear system with $R$ and $R^T$ is, 
\[
t_{\text{QR, apply}} \in  \mathcal{O}\Bigg(\sum_{l=1}^{L}2^l\Big(2^{-2l/3}N^{2/3}\Big)^2\Bigg) = \mathcal{O}\big(N^{4/3}\big)
\]

\paragraph{spaQR} Let us now estimate the complexity of the spaQR factorization, by making additional assumptions on the sparsification process. 
\begin{enumerate}
    \item[5.] Sparsification reduces the size of an interface at level $l$ to, 
\[c_l' \in \mathcal{O}(2^{-l/3}N^{1/3})\]
This means that the rank scales roughly as the diameter of the separator. 
\item[6.] An interface has $\mathcal{O}(1)$ neighbor interfaces
\item[7.] The aspect ratio of the diagonal blocks corresponding to each interface is $\Theta(\alpha)$. This implies that the number of rows of an interface at level $l$ is, 
\[
r_l' \in \mathcal{O}(\alpha 2^{-l/3}N^{1/3})
\]

\end{enumerate}
Assumption 5 is a consequence of low rank interactions between separators that are far away in $G_{A^TA}$. This is comparable to complexity assumptions in the fast multipole method~\cite{FMM_1,greengard_rokhlin_1997}, spaND~\cite{2019arXiv190102971C}, and HIF~\cite{Ho2016HierarchicalIF}. Assumption 7 is a consequence of the extra sparsification step defined in \Cref{sec: sparsify}. 

The fill-in in the sparsified QR process results in at most $\mathcal{O}(2^{-l/3}N^{1/3})$ entries in each row and column. This is in part due to the assumption on the size of the interfaces, the number of neighbor interfaces and the fact that new connections are only made between distance 1 neighbors of a node in $G_{A^TA}$. 

The cost of the spaQR factorization is split into four parts:
\begin{itemize}
    \item \textit{Householder QR on interiors and separators.} A separator block has $c_l'$ columns and $r_l'$ rows right before it is factorized. Further, it has at most $\mathcal{O}(2^{-l/3}N^{1/3})$ non-zeros per row/column and $\mathcal{O}(1)$ neighbor interfaces. Then, the cost of doing QR on a separator is \[
    h_l' \in \mathcal{O}\big(\alpha (2^{-l/3}N^{1/3})^3 \big) = \mathcal{O}\big(2^{-l}N\big)
    \]
    \item \textit{Scaling of interfaces.} The cost of scaling (QR on a block of size $r_l'\times c_l'$) an interface is $\mathcal{O}\big(\alpha2^{-l}N\big)$.
    \item \textit{Sparsification of interfaces: Step 1} The rows below the diagonal of each remaining interface are sparsified in the first step of the sparsification stage. In this step, we perform a rank-revealing QR on a block with $\mathcal{O}\big((\alpha-1)c_l')$ rows and $\mathcal{O}(n_{\text{nbr}}c_l')$  columns, where $n_{\text{nbr}}$ is the number of neighbor interfaces. We have assumed that the number of neighbor interfaces are $\mathcal{O}(1)$. Hence, for large $\alpha$, the cost of sparsifying the rows below the diagonal of an interface is, 
    \[
        e_l' \in \mathcal{O}\big( (\alpha-1) (2^{-l/3}N^{1/3})^3 \big) = \mathcal{O}(\alpha 2^{-l} N)
    \]
    \item \textit{Sparsification of interfaces: Step 2} In this step, we compute a  rank-revealing QR on the block $\begin{bmatrix}
    A_{np}^T & A_{p_1n}
    \end{bmatrix}$. The cost of rank-revealing QR is $\mathcal{O}(mnr)$ where $r \leq \min(m,n)$. The number of rows of this block is $c_l'$. The number of columns is $\mathcal{O}(\alpha 2^{-l/3}N^{1/3})$ due to our assumptions on the size of the interface and the number of neighbor interfaces. Thus, this step costs, 
    \[
        s_l' \in \mathcal{O}\big(\alpha (2^{-l/3}N^{1/3})^3 \big) = \mathcal{O}(\alpha 2^{-l} N)
    \] 
    
\end{itemize}
Hence, the total cost of the spaQR algorithm is, 
\begin{align*}
    t_{\text{spaQR}} \in \mathcal{O}\Bigg(\sum_{l=1}^{L}\alpha 2^{l}2^{-l}N\Bigg) &= \mathcal{O}\Bigg(\sum_{l=1}^{L} \alpha N\Bigg) \\ &= \mathcal{O}(\alpha N \log N), \qquad L \in \Theta(\log(N/N_0)) \\
    &=\mathcal{O}(M \log N)
\end{align*}

The cost of applying the factorization can be derived similarly to the analysis for the multifrontal QR algorithm. Solving the least squares system consists essentially in solving $W^T W x = A^T b$, where $W$ is a sequence of sparse orthogonal and upper-triangular matrices. Using our assumption on the size of the interfaces, the total cost of applying the spaQR factorization is, 
\[
t_{\text{spaQR, apply}}  \in \mathcal{O}\Bigg(\sum_{l=1}^{L}\alpha2^{l} \Big(2^{-l/3}N^{1/3}\Big)^2 \Bigg) = \mathcal{O}(\alpha N) = \mathcal{O}(M)
\]
The memory required to store the factorization scales as the cost of applying the factorization. Hence, 
\[
\text{mem}_{\text{spaQR}} = \mathcal{O}(M)
\]
We show numerical results in \Cref{sec: profiling} that corroborate the assumptions made here. 

\section{Benchmarks}
\label{sec: benchmarks}

In this section, we benchmark the performance of the algorithm for solving sparse linear least squares problem arising in various applications. The spaQR algorithm is used to compute an approximate factorization using a user-defined tolerance $\epsilon$. The approximate factorization is then used as a preconditioner with the CGLS iterative solver. The convergence criteria is set as $\|A^T(Ax-b)\|_{_2}/\|A^Tb\|_{_2}\leq 10^{-12}$. The performance of spaQR is compared with direct multifrontal QR (computed using the spaQR code with no compression but otherwise the same parameters), which is denoted by `Direct'. We also compare against a CGLS preconditioned using a standard diagonal preconditioner (using the Eigen library's rountine), which is denoted by `Diag.'. 

The code was written in C++. We use GCC 8.1.0 and Intel(R) MKL 2019 for Linux for the BLAS and LAPACK operations. The number of levels in the nested dissection process is chosen as $\lceil\log (N/64)/\log 2\rceil$ for a matrix of size $M \times N$. Low rank approximations are performed using LAPACK's dlaqps routine which performs a column pivoted QR on $r$ columns. The value $r$ is chosen such that $\frac{|R_{ii}|}{|R_{11}|} \geq \epsilon$ for $1\leq i \leq r$, where $R$ is the upper triangular matrix that comes out of the column pivoted QR method. We typically begin sparsification on levels 3 or 4.

\subsection{Scaling with problem size}
\label{subsec:scaling_w_size}

First, we study the performance of the algorithm with increasing problem size on the Inverse Poisson problem defined on uniform 2D and 3D grids. We use a regular geometric partitioning on $A^TA$ to get the separators and infer the interfaces as outlined in \Cref{ord_clus}.

The PDE-constrained optimization problem is as follows,
\begin{align*}
    \min_{u,z} \quad & \frac{1}{2}\int_{\Omega} (u-u_d)^2 dx +\frac{1}{2} \lambda \int_{\Omega}z^2dx\\
    \text{subject to} \quad & -\nabla \cdot (z \nabla u) = h \text{ in }\Omega\\
     \quad & u = 0 \text{ in } d\Omega
\end{align*}

The variable coefficient Poisson equation is discretized using a finite difference scheme on a staggered grid~\cite{FDM_invpoi}. We want to recover the variables $u$ and the diffusion coefficients $z$ at every grid point, using the available observed states $u_d$. One approach to solving this optimization problem~\cite{estrin2019implementing} involves solving the least squares problem, 
\[J^T x = b\]
where $J$ is the Jacobian of the discretized constraints. The Jacobian can be written as
$J = \begin{bmatrix}
        J_u & J_z
\end{bmatrix}
$ where the columns of $J_u$ (resp.\ $J_z$) correspond to a partial derivative with respect to one of the $u_i$ ($z_i$) variables. 

The Jacobian $J$ can be evaluated with different values of $(u,z)$. Generally speaking, if we choose random values for $u$ and $z$ we will get a least squares matrix ($J^T$) with an aspect ratio of 2. It is possible to reduce the aspect ratio of the matrix by artificially creating rows that are equal to 0. Then, by removing zero rows, we can reduce the aspect ratio of matrix $J^T$. A detailed explanation of the process of generating matrices with different aspect ratio is given in \Cref{sec: mat_gen_2d,sec: mat_gen_3d}.

\subsubsection{2D Inverse Poisson problem}

For a 2D staggered grid, the number of columns in the least squares problem is $N=n^2$. The aspect ratio of the matrix can be varied roughly from 1 to 2 as described in \Cref{sec: mat_gen_2d}. In our benchmarks, we consider least squares matrix $J^T$ with $N=n^2$ columns and aspect ratios $\alpha \approx 2$, 1.5, and 1.05.
\begin{figure}[thbp]
\centering
\begin{subfigure}{\textwidth}
\scalebox{0.7}{
\begin{tikzpicture}
\pgfplotsset{
every axis legend/.append style={
at={(0.5,1.03)},
anchor=south
},
}
\begin{loglogaxis}[
legend columns = 2,
legend style = {draw = none},
    scale = 0.65,
    ylabel={Time to factorize (s)},
    ymin=0.15, ymax=5000,
    xtick = \empty,
    extra x ticks = {16000, 6.5*10^4, 2.5*10^5, 10^6, 4*10^6},
   extra x tick labels = \empty,
    ymajorgrids=true,
    line width=0.25mm,
    grid style=dashed,
    legend entries = { $\epsilon = 10^{-4}$, Direct}
]
\addplot coordinates {
    (128^2, 0.67)(256^2,1.2) (512^2, 3.5) (1024^2, 13.5) (2048^2, 55)
    };
\addplot+ [black,mark = triangle*, mark options ={fill = black}] coordinates {
    (128^2, 0.26)(256^2,1) (512^2, 9.67) (1024^2, 43.3) (2048^2, 292)
    };

\addplot [black, domain = 128^2:2048^2] {x/150000};
\node [ anchor=center] at (100*10^4,2) {$\mathcal{O}(N)$};
\addplot [black, dashed, domain = 128^2:2048^2] {x^(3/2)/10000000};
\node [ anchor=center] at (0.2*10^6,100) {$\mathcal{O}(N^{3/2})$};
\end{loglogaxis}
\end{tikzpicture}}
\scalebox{0.7}{
\begin{tikzpicture}
\pgfplotsset{
every axis legend/.append style={
at={(0.5,1.03)},
anchor=south
},
}
\begin{loglogaxis}[
legend columns = 2,
legend style = {draw = none},
    scale = 0.65,
    ymin=0.15, ymax= 5000,
    xtick = \empty,
    extra x ticks = {16000, 6.5*10^4, 2.5*10^5, 10^6, 4*10^6},
    extra x tick labels = \empty,
    yticklabel = \empty,
    line width=0.25mm,
    ymajorgrids=true,
    grid style=dashed,
    legend entries = { $\epsilon = 10^{-2}$, Direct}
]
\addplot coordinates {
    (128^2, 0.6)(256^2,1.6) (512^2, 7) (1024^2, 33) (2048^2, 149)
    };
\addplot+ [black,mark = triangle* , mark options ={fill = black}] coordinates {
    (128^2, 0.5)(256^2,9) (512^2, 114.7) (1024^2, 1838) (2048^2, nan)
    };
\addplot+ [draw=black, dotted, mark= triangle, mark options={solid}] coordinates {
    (128^2, 0.5)(256^2,9) (512^2, 114.7) (1024^2, 1838) (2048^2, 29408)
};
\addplot [black, domain = 128^2:2048^2] {x/70000};
\node [ anchor=center] at (100*10^4,4) {$\mathcal{O}(N)$};
\addplot [black, dashed, domain = 128^2:2048^2] {x^(2)/(2*10^8)};
\node [ anchor=center] at (0.05*10^6,200) {$\mathcal{O}(N^{2})$};
\end{loglogaxis}
\end{tikzpicture}}
\scalebox{0.7}{
\begin{tikzpicture}
\pgfplotsset{
every axis legend/.append style={
at={(0.5,1.03)},
anchor=south
},
}
\begin{loglogaxis}[
legend columns = 2,
legend style = {draw = none},
    scale = 0.65,
    ymin=0.15, ymax=5000,
    xtick = \empty,
    extra x ticks = {16000, 6.5*10^4, 2.5*10^5, 10^6, 4*10^6},
    extra x tick labels = {16k, , 0.25M, ,4M},
    yticklabel = \empty,
    extra x tick labels = \empty,
    line width=0.25mm,
    ymajorgrids=true,
    grid style=dashed,
    legend entries = { $\epsilon = 10^{-2}$, Direct}
]
\addplot coordinates {
    (128^2,0.9)(256^2,2.7) (512^2, 11.9) (1024^2, 54) (2048^2,250)
    };
\addplot+ [black, mark = triangle*, mark options ={fill = black}] coordinates {
    (128^2, 1.4)(256^2,16.6) (512^2, 221) (1024^2, 3305) (2048^2, nan)
    };
\addplot+ [draw=black, dotted, mark= triangle, mark options={solid}] coordinates {
    (128^2, 1.4)(256^2,16.6) (512^2, 221) (1024^2, 3305) (2048^2, 52880)
};
\addplot [black, domain = 128^2:2048^2] {x/40000};
\node [ anchor=center] at (100*10^4,9) {$\mathcal{O}(N)$};
\addplot [black, dashed, domain = 128^2:2048^2] {x^(2)/(1*10^8)};
\node [ anchor=center] at (0.05*10^6,200) {$\mathcal{O}(N^{2})$};
\end{loglogaxis}
\end{tikzpicture}}
\end{subfigure}

\begin{subfigure}{\textwidth}
\scalebox{0.7}{
\begin{tikzpicture}
\begin{loglogaxis}[
    scale = 0.65,
    ylabel={$\#$ CGLS},
    ymin=1, ymax=100,
    xtick = \empty,
    extra x ticks = {16000, 6.5*10^4, 2.5*10^5, 10^6, 4*10^6},
     extra x tick labels = {16k, , 0.25M, ,4M},
    xlabel={N, $\alpha \approx 1.05$},
    ymajorgrids=true,
    line width=0.25mm,
    grid style=dashed
]

\addplot coordinates {
    (128^2, 3)(256^2,4) (512^2, 5) (1024^2, 8) (2048^2, 14)
    };

\end{loglogaxis}
\end{tikzpicture}}
\scalebox{0.7}{
\begin{tikzpicture}
\begin{loglogaxis}[
    scale = 0.65,
    ymin=1, ymax=100,
    xtick = \empty,
    extra x ticks = {16000, 6.5*10^4, 2.5*10^5, 10^6, 4*10^6},
     extra x tick labels = {16k, , 0.25M, ,4M},
    xlabel={N, $\alpha \approx 1.5$},
    yticklabel = \empty,
    ymajorgrids=true,
    line width=0.25mm,
    grid style=dashed
]

\addplot coordinates {
    (128^2, 10)(256^2,10) (512^2, 14) (1024^2, 21) (2048^2, 38)
    };

\end{loglogaxis}
\end{tikzpicture}}
\scalebox{0.7}{
\begin{tikzpicture}
\begin{loglogaxis}[
    scale = 0.65,
    ymin=1, ymax=100,
    xtick = \empty,
    extra x ticks = {16000, 6.5*10^4, 2.5*10^5, 10^6, 4*10^6},
     extra x tick labels = {16k, , 0.25M, ,4M},
    xlabel={N, $\alpha \approx 2$},
    yticklabel = \empty,
    ymajorgrids=true,
    line width=0.25mm,
    grid style=dashed
]

\addplot coordinates {
    (128^2, 7)(256^2,8) (512^2, 10) (1024^2, 14) (2048^2, 20)
    };
\end{loglogaxis}
\end{tikzpicture}}
\end{subfigure}
\caption{Results for a 2D Inverse Poisson least squares problem for varying values of the aspect ratio $\alpha=M/N$. For small enough values of $\epsilon$, the iteration counts increase slowly.  Note that for smaller aspect ratio, the condition number is higher and hence requires a higher accuracy (smaller $\epsilon$) to converge in $< 50$ iterations. The factorization time scales as $\mathcal{O}(N)$ for all three values of $\alpha$. }
\label{2d_invpoi}
\end{figure}
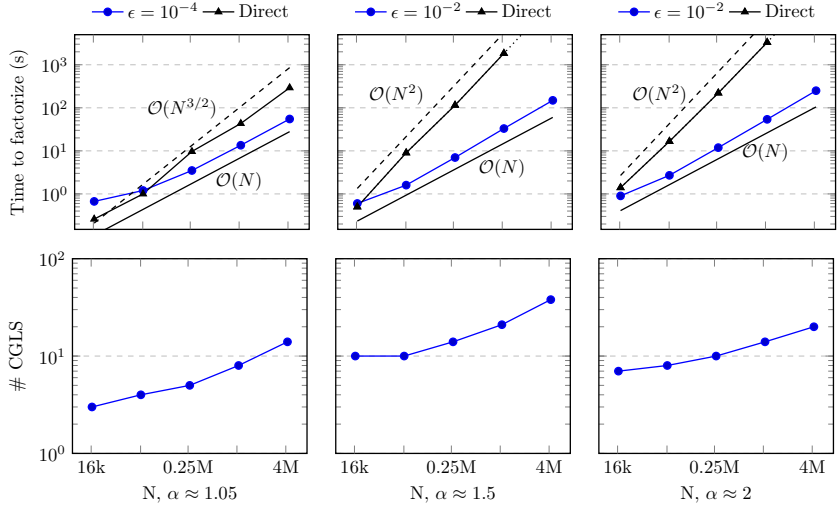
\begin{figure}[thbp]
\centering
\begin{subfigure}{\textwidth}
\centering
\scalebox{0.7}{
\begin{tikzpicture}
\pgfplotsset{
every axis legend/.append style={
at={(0.5,1.03)},
anchor=south
},
}
\begin{semilogyaxis}[
legend columns = 2,
legend style = {draw = none},
legend cell align={left},
    scale = 0.65,
    ylabel={Residual},
    xlabel={Iterations, $\alpha\approx1.05$},
    ymin = 1e-14, ymax =1,
    ymajorgrids=true,
    line width=0.25mm,
    grid style=dashed,
    legend entries = {spaQR  $\epsilon = 10^{-4}$, Diag.}
]

\addplot table {results/2d_invpoi_1_spaQR.dat};
\addplot table {results/2d_invpoi_1_cgls.dat};
\end{semilogyaxis}
\end{tikzpicture}}
\scalebox{0.7}{
\begin{tikzpicture}
    \pgfplotsset{
every axis legend/.append style={
at={(0.5 ,1.03)},
anchor=south
},
}
\begin{semilogyaxis}[
legend columns = 2,
legend style = {draw = none},
legend cell align={left},
    scale = 0.65,
    xlabel={Iterations, $\alpha\approx1.5$},
    ymin = 1e-14, ymax =1,
    yticklabel = \empty,
    ymajorgrids=true,
    line width=0.25mm,
    grid style=dashed,
    legend entries = { spaQR $\epsilon = 10^{-2}$, Diag.}
]

\addplot table {results/2d_invpoi_1_5_spaQR.dat};
\addplot table {results/2d_invpoi_1_5_cgls.dat};
\end{semilogyaxis}
\end{tikzpicture}
}
\scalebox{0.7}{
\begin{tikzpicture}
\pgfplotsset{
every axis legend/.append style={
at={(0.5 ,1.03)},
anchor=south
},
}
\begin{semilogyaxis}[
legend columns = 2,
legend style = {draw = none},
legend cell align={left},
    scale = 0.65,
    xlabel={Iterations, $\alpha\approx 2$},
    ymin = 1e-14, ymax =1,
    yticklabel = \empty,
    ymajorgrids=true,
    line width=0.25mm,
    grid style=dashed,
    legend entries = { spaQR  $\epsilon = 10^{-2}$, Diag.}
]

\addplot table {results/2d_invpoi_2_spaQR.dat};
\addplot table {results/2d_invpoi_2_cgls.dat};
\end{semilogyaxis}
\end{tikzpicture}
}
\end{subfigure}
\caption{Convergence of the residual $\|A^{T}(Ax-b)\|_{_2}/\|A^Tb\|_{_2}$ with the number of CGLS iterations with the diagonal preconditioner and spaQR for the $2048 \times 2048$ 2D Inverse Poisson problem. Note the superiority of CGLS preconditioned with the spaQR algorithm for all three values of the aspect ratio $\alpha$.}
\label{fig: invpoi_cgls}
\end{figure}
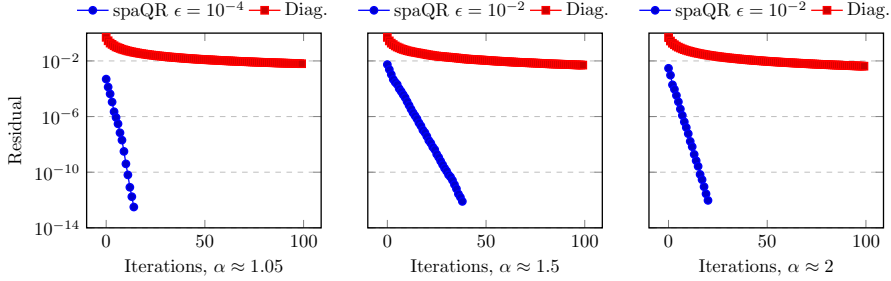
The results for the 2D Inverse Poisson problem are shown in \Cref{2d_invpoi}; the number of columns are $N = n^2$ and the number of rows are $M = \alpha N$, with the aspect ratio $\alpha$ varying from $1.05$ to $2$. For small enough $\epsilon$, the number of CGLS iterations needed to converge increase slowly with increasing problem size. We note that for problems with $\alpha \approx 1.05$, we need a smaller value of $\epsilon$ ($10^{-4}$) for good convergence, as these problems are ill-conditioned.

In the case of direct multifrontal QR (with the same nested dissection ordering), the size of the top separator will be $\mathcal{O}(N^{1/2})$. As more and more fill-in are introduced as the factorization proceeds, the block corresponding to the top separator will be dense when we reach the top of the elimination tree. The size of the top separator block will be of size $\mathcal{O}\big((\alpha-1)N+ N^{1/2}\big) \times \mathcal{O}(N^{1/2})$ and a Householder QR on this block will cost $\mathcal{O}\big((\alpha-1)N^2+ N^{3/2}\big)$ flops. On the other hand, the factorization time in spaQR scales as $\mathcal{O}(N)$ for all three values of $\alpha$ as expected (see \Cref{2d_invpoi}). 

We also compare the convergence of the CGLS residual with a diagonal preconditioner and using spaQR as the preconditioner. In \Cref{fig: invpoi_cgls}, we show the convergence results for the least squares problem on a $2048\times 2048$ grid for three values of the aspect ratio $\alpha$. In all three cases, the spaQR preconditioned CGLS takes less than 30 iterations to converge to a residual of $10^{-12}$. However, the diagonal preconditioner performs poorly, taking around 100 iterations to reach a residual of $10^{-2}$.

\subsubsection{3D Inverse Poisson problem}
\label{subsec: 3d_invpoi}

For a 3D $n \times n \times n$ grid, the least squares problem has $N = n^3$ columns and the aspect ratio can be varied from 1 to 2 as described in \Cref{sec: mat_gen_3d}. In our benchmarks, we consider problems with aspect ratio $\alpha \in \{1.05, 1.5, 2\}$.

\begin{figure}[thbp]
\centering
\begin{subfigure}{\textwidth}
\centering
\scalebox{0.7}{
\begin{tikzpicture}
\pgfplotsset{
every axis legend/.append style={
at={(0.5,1.03)},
anchor=south
},
}
\begin{loglogaxis}[
legend columns = 2,
legend style = {draw = none},
    scale = 0.65,
    ylabel={Time to factorize (s)},
    ymin=3, ymax=20000,
    xtick = \empty,
    xtick = {3*10^4, 1*10^5, 3*10^5, 5*10^5, 8*10^5, 2*10^6},
    xticklabels = \empty,
    ymajorgrids=true,
    line width=0.25mm,
    grid style=dashed,
    legend entries = { $\epsilon = 10^{-2}$, Direct}
]
\addplot coordinates {
    (32^3, 6.2) (48^3,35) (64^3,124) (80^3, 358) (96^3, 752) (128^3, 2542)
    };
\addplot+ [black,mark = triangle*, mark options={fill=black}] coordinates {
    (32^3, 14.3) (48^3,150) (64^3,829) (80^3, 3581) (96^3, nan) (128^3, nan)
    };
\addplot+ [draw=black,  dotted, mark= triangle, mark options={solid}] coordinates {
    (32^3, 14.3) (48^3,150) (64^3, 829) (80^3, 3581) (96^3, 10692) (128^3, 60000)
};
\addplot [black, dashed, domain = 30000:2*10^6] {x^(2)/(4*10^7)};
\node [anchor=center] at (20*10^4,5000) {$\mathcal{O}(N^{2})$};

\addplot [black, domain = 30000:2*10^6] {x^(1.4)/600000};
\node [anchor=center] at (1*10^6,100) {$\mathcal{O}(N^{1.4})$};
\end{loglogaxis}
\end{tikzpicture}}
\scalebox{0.7}{
\begin{tikzpicture}
\pgfplotsset{
every axis legend/.append style={
at={(0.5,1.03)},
anchor=south
},
}
\begin{loglogaxis}[
legend columns = 2,
legend style = {draw = none},
    scale = 0.65,
    ymin=3, ymax=20000,
    xtick = {3*10^4, 1*10^5, 3*10^5, 5*10^5, 8*10^5, 2*10^6},
    xticklabels = \empty,
    yticklabel = \empty,
    line width=0.25mm,
    ymajorgrids=true,
    grid style=dashed,
    legend entries = { $\epsilon = 10^{-2}$, Direct}
]
\addplot coordinates {
    (32^3, 8.7) (48^3,56) (64^3,216) (80^3, 594) (96^3, 1169) (128^3, 4124)
    };
\addplot+ [black,mark = triangle*, mark options={fill=black}] coordinates {
    (32^3, 31.8) (48^3,470) (64^3, 3453) (80^3, nan) (96^3, nan) (128^3, nan)
    };
\addplot+ [draw=black,  dotted, mark= triangle, mark options={solid}] coordinates {
    (32^3, 31.8) (48^3,470) (64^3, 3453) (80^3, 16465) (96^3, nan) (128^3, nan)
};
\addplot [black, dashed, domain = 30000:2*10^6] {x^(7/3)/(5*10^8)};
\node [anchor=center] at (5*10^4,5000) {$\mathcal{O}(N^{7/3})$};
\addplot [black, domain = 30000:2*10^6] {x^(1.5)/1000000};
\node [anchor=center] at (0.9*10^6,200) {$\mathcal{O}(N^{1.5})$};
\end{loglogaxis}
\end{tikzpicture}}
\scalebox{0.7}{
\begin{tikzpicture}
\pgfplotsset{
every axis legend/.append style={
at={(0.5,1.03)},
anchor=south
},
}
\begin{loglogaxis}[
legend columns = 2,
legend style = {draw = none},
    scale = 0.65,
    ymin=3, ymax=20000,
    xtick = {3*10^4, 1*10^5, 3*10^5, 5*10^5, 8*10^5, 2*10^6},
    xticklabels = \empty,
    yticklabel = \empty,
    line width=0.25mm,
    ymajorgrids=true,
    grid style=dashed,
    legend entries = { $\epsilon = 10^{-2}$, Direct}
]
\addplot coordinates {
    (32^3, 15) (48^3,93) (64^3,372) (80^3,941) (96^3, 2118) (128^3, 9093)
    };
\addplot+ [black,mark = triangle*, mark options={fill=black}] coordinates {
    (32^3, 66) (48^3,1008) (64^3,8840) (80^3, nan) (96^3, nan) (128^3, nan)
    };
\addplot+ [draw=black, dotted, mark= triangle, mark options={solid}] coordinates {
    (32^3, 66) (48^3,1008) (64^3, 8840) (80^3, 42000) (96^3, nan) (128^3, nan)
};
\addplot [black, dashed, domain = 30000:2*10^6] {x^(7/3)/(2*10^8)};
\node [anchor=center] at (5*10^4,5000) {$\mathcal{O}(N^{7/3})$};

\addplot [black, domain = 30000:2*10^6] {x^(1.5)/600000};
\node [anchor=center] at (0.9*10^6,300) {$\mathcal{O}(N^{1.5})$};
\end{loglogaxis}
\end{tikzpicture}}
\end{subfigure}

\begin{subfigure}{\textwidth}
\centering
\scalebox{0.7}{
\begin{tikzpicture}
\begin{loglogaxis}[
    scale = 0.65,
    xlabel={N, $\alpha \approx 1.05$},
    ylabel={$\#$ CGLS},
    xtick = \empty,
    extra x ticks = {3*10^4, 1*10^5, 3*10^5, 5*10^5, 8*10^5, 2*10^6},
    extra x tick labels = {30k, ,0.3M, ,,2M},
    ymin=1, ymax = 100,
    ymajorgrids=true,
    line width=0.25mm,
    grid style=dashed
]
\addplot coordinates {
    (32^3, 10) (48^3,12) (64^3,15) (80^3, 22) (96^3, 25) (128^3, 41)
    };
\end{loglogaxis}
\end{tikzpicture}}
\scalebox{0.7}{
\begin{tikzpicture}
\begin{loglogaxis}[
    scale = 0.65,
    ymin=1,ymax=100,
    xlabel={N, $\alpha \approx 1.5$},
    xtick = \empty,
    extra x ticks = {3*10^4, 1*10^5, 3*10^5, 5*10^5, 8*10^5, 2*10^6},
   extra x tick labels = {30k, ,0.3M, ,,2M},
    yticklabel = \empty,
    ymajorgrids=true,
    line width=0.25mm,
    grid style=dashed
]
\addplot coordinates {
    (32^3, 8) (48^3,8) (64^3,9) (80^3, 9) (96^3, 9) (128^3, 9)
    };
\end{loglogaxis}
\end{tikzpicture}}
\scalebox{0.7}{
\begin{tikzpicture}
\begin{loglogaxis}[
    scale = 0.65,
    ymin=1, ymax=100,
    xtick = \empty,
    xlabel={N, $\alpha \approx 2$},
    extra x ticks = {3*10^4, 1*10^5, 3*10^5, 5*10^5, 8*10^5, 2*10^6},
    extra x tick labels = {30k, ,0.3M, ,,2M},
    yticklabel = \empty,
    ymajorgrids=true,
    line width=0.25mm,
    grid style=dashed
]
\addplot coordinates {
    (32^3, 8) (48^3,10) (64^3,10) (80^3, 11) (96^3, 11) (128^3, 12)
    };
\end{loglogaxis}
\end{tikzpicture}}
\end{subfigure}
\caption{Results of 3D Inverse Poisson least squares problem for varying values of aspect ratio,  $\alpha=M/N$. For small enough values of $\epsilon$, the iteration counts increase slowly. Empirically, the number of flops roughly scales between $\mathcal{O}(N^{1.4})$ to $\mathcal{O}(N^{1.5})$ as shown. The direct multifrontal QR scales like $\mathcal{O}\big((\alpha-1)N^{7/3}+N^2\big)$ as expected.}
\label{3d_invpoi}
\end{figure}
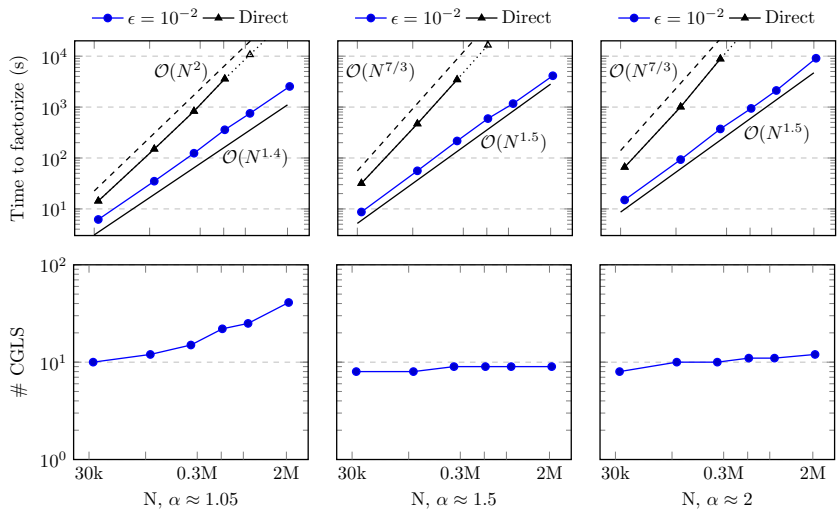

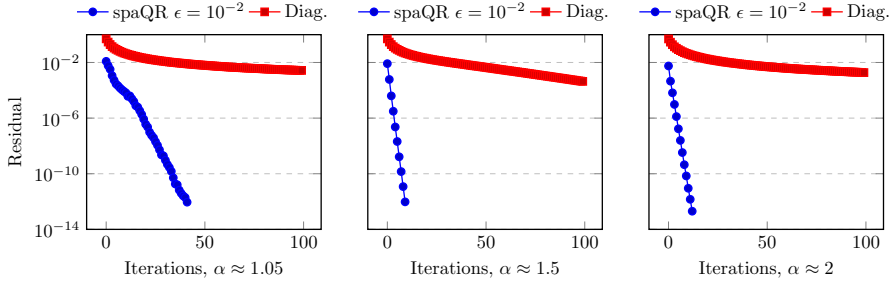
\begin{figure}[htbp]
\centering
\begin{subfigure}{\textwidth}
\centering
\scalebox{0.7}{
\begin{tikzpicture}
\pgfplotsset{
every axis legend/.append style={
at={(0.5,1.03)},
anchor=south
},
}
\begin{semilogyaxis}[
legend columns = 2,
legend style = {draw = none},
legend cell align={left},
    scale = 0.65,
    ylabel={Residual},
    xlabel={Iterations, $\alpha\approx1.05$},
    ymin = 1e-14, ymax =1,
    ymajorgrids=true,
    line width=0.25mm,
    grid style=dashed,
    legend entries = {spaQR  $\epsilon = 10^{-2}$, Diag.}
]

\addplot table {results/3d_invpoi_1_spaQR.dat};
\addplot table {results/3d_invpoi_1_cgls.dat};
\end{semilogyaxis}
\end{tikzpicture}}
\scalebox{0.7}{
\begin{tikzpicture}
    \pgfplotsset{
every axis legend/.append style={
at={(0.5 ,1.03)},
anchor=south
},
}
\begin{semilogyaxis}[
legend columns = 2,
legend style = {draw = none},
legend cell align={left},
    scale = 0.65,
    xlabel={Iterations, $\alpha\approx1.5$},
    ymin = 1e-14, ymax =1,
    yticklabel = \empty,
    ymajorgrids=true,
    line width=0.25mm,
    grid style=dashed,
    legend entries = {spaQR $\epsilon = 10^{-2}$, Diag.}
]

\addplot table {results/3d_invpoi_1.5_spaQR.dat};
\addplot table {results/3d_invpoi_1.5_cgls.dat};
\end{semilogyaxis}
\end{tikzpicture}
}
\scalebox{0.7}{
\begin{tikzpicture}
\pgfplotsset{
every axis legend/.append style={
at={(0.5 ,1.03)},
anchor=south
},
}
\begin{semilogyaxis}[
legend columns = 2,
legend style = {draw = none},
legend cell align={left},
    scale = 0.65,
    xlabel={Iterations, $\alpha\approx 2$},
    ymin = 1e-14, ymax =1,
    yticklabel = \empty,
    ymajorgrids=true,
    line width=0.25mm,
    grid style=dashed,
    legend entries = { spaQR  $\epsilon = 10^{-2}$, Diag.}
]

\addplot table {results/3d_invpoi_2_spaQR.dat};
\addplot table {results/3d_invpoi_2_cgls.dat};
\end{semilogyaxis}
\end{tikzpicture}
}
\end{subfigure}
\caption{Convergence of the residual $\|A^{T}(Ax-b)\|_{_2}/\|A^Tb\|_{_2}$ with the number of CGLS iterations with spaQR and a diagonal preconditioner for the $128 \times 128 \times 128$ 3D Inverse Poisson least squares problem. Note the superiority of CGLS preconditioned with the spaQR method for all three values of the aspect ratio $\alpha$. }
\label{fig: 3d_invpoi_cgls}
\end{figure}

\Cref{3d_invpoi} shows the number of CGLS iterations to converge and the time of factorization for the Inverse Poisson problem on a 3D grid. Theoretically, we expect the factorization time to scale as $\mathcal{O}(N \log N)$ for a fixed value of $\alpha$ (see \Cref{sec:complexity}). However, the empirical cost of factorization scales as $\mathcal{O}(N^{1.4})$ for $\alpha \approx 1.05$ and as $\mathcal{O}(N^{1.5})$ for $\alpha \approx$ 1.5, 2. Based on other benchmarks, we believe this is due to non-asymptotic effects (e.g., terms in the computational runtime that vanish asymptotically but may dominate at small sizes). The number of iterations increases slowly (almost a constant for $\alpha=1.5, 2$) for the range of problems considered. In \Cref{fig: 3d_invpoi_cgls} we compare the performance of CGLS with a standard diagonal preconditioner and spaQR. These results show the superior performance of spaQR as compared to some of the standard techniques in this area.

\subsection{Profiling}
\label{sec: profiling}

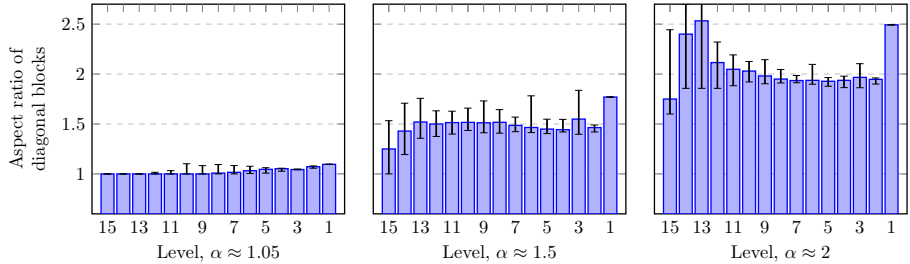
\begin{figure}
    \centering
\scalebox{0.7}{
    \begin{tikzpicture}
        \begin{axis}[
        ylabel style = {align=center},
            scale = 0.7,
            ylabel = {Aspect ratio of \\ diagonal blocks},
            xlabel = {Level, $\alpha\approx 1.05 $},
            ymin = 0.6, ymax = 2.7,
            x dir=reverse,
             xticklabel = \empty,
            extra x ticks = { 15,14,13,12, 11, 10, 9,8,7,6,5,4,3,2,1 },
            extra x tick labels = {15,,13,,11,,9,,7,,5,,3,,1},
            enlarge x limits=0.07,
            bar width = 7pt,
            ymajorgrids=true,
            line width =0.25mm,
            grid style=dashed,
            xtick align = inside
        ]
        \addplot[ybar, draw=blue, fill=blue!30,error bars/.cd,
                    y explicit,
                    y dir=both, 
                    error bar style={line width=0.3mm, black}] 
                    table [
                    x = level,
                    y = median,
                    y error plus expr=\thisrow{q2}-\thisrow{median},
                    y error minus expr=\thisrow{median}-\thisrow{q1}
                    ] 
                    {results/3d_invpoi_1_alpha_level.dat};
        \end{axis}
    \end{tikzpicture}}
\scalebox{0.7}{
    \begin{tikzpicture}
        \begin{axis}[
            scale = 0.7,
            ymin = 0.6, ymax = 2.7,
            xlabel = {Level, $\alpha\approx 1.5 $},
            x dir=reverse,
             xticklabel = \empty,
            extra x ticks = { 15,14,13,12, 11, 10, 9,8,7,6,5,4,3,2,1 },
            extra x tick labels = {15,,13,,11,,9,,7,,5,,3,,1},
            yticklabel = \empty,
            enlarge x limits=0.07,
            bar width = 7pt,
            ymajorgrids=true,
            line width =0.25mm,
            grid style=dashed,
            xtick align = inside
        ]
        \addplot[ybar, draw=blue, fill=blue!30,error bars/.cd,
                    y explicit,
                    y dir=both, 
                    error bar style={line width=0.3mm, black}] 
                    table [
                    x = level,
                    y = median,
                    y error plus expr=\thisrow{q2}-\thisrow{median},
                    y error minus expr=\thisrow{median}-\thisrow{q1}
                    ] 
                    {results/3d_invpoi_1.5_alpha_level.dat};
        \end{axis}
    \end{tikzpicture}}
\scalebox{0.7}{
    \begin{tikzpicture}
        \begin{axis}[
            scale = 0.7,
            ymin = 0.6, ymax = 2.7,
            xlabel = {Level, $\alpha\approx 2 $},
            x dir=reverse,
            xticklabel = \empty,
            extra x ticks = { 15,14,13,12, 11, 10, 9,8,7,6,5,4,3,2,1 },
            extra x tick labels = {15,,13,,11,,9,,7,,5,,3,,1},
            yticklabel = \empty,
            enlarge x limits=0.07,
            bar width = 7pt,
            ymajorgrids=true,
            line width =0.25mm,
            grid style=dashed,
            xtick align = inside
        ]
        \addplot[ybar, draw=blue, fill=blue!30,error bars/.cd,
                    y explicit,
                    y dir=both, 
                    error bar style={line width=0.3mm, black}] 
                    table [
                    x = level,
                    y = median,
                    y error plus expr=\thisrow{q2}-\thisrow{median},
                    y error minus expr=\thisrow{median}-\thisrow{q1}
                    ] 
                    {results/3d_invpoi_2_alpha_level.dat};
        \end{axis}
    \end{tikzpicture}}
\caption{The median aspect ratio of the diagonal blocks corresponding to the interfaces per level for the 3D Inverse Poisson problem on $N=128^3$. The error bars show the interquartile range.}
    \label{fig:3d_aspect_ratio}
\end{figure}

In this section, we try to understand the time and memory requirements of the spaQR algorithm. First, we show some experimental evidence to back up the assumptions made in the complexity analysis. Assumptions 1 to 6 are the same as the assumptions used to get the complexity estimate for squares matrices in~\cite{gnanasekaran2020hierarchical} in which we also provide experimental evidence for them. \Cref{fig:3d_aspect_ratio} shows that the aspect ratio of the diagonal blocks is almost a constant at each level in the elimination tree and is bounded by $\Theta(\alpha)$ at every level. 

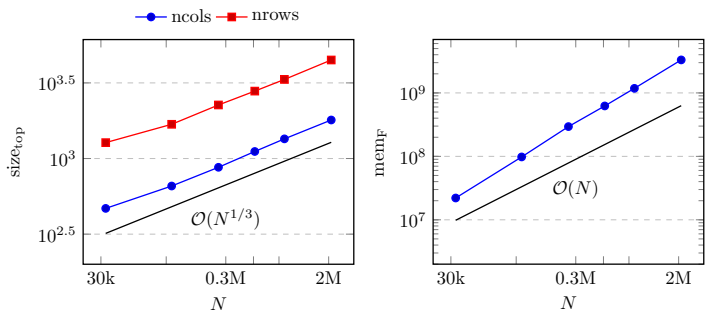
\begin{figure}
    \centering
    \scalebox{0.7}{
\begin{tikzpicture}
\pgfplotsset{
every axis legend/.append style={
at={(0.5 ,1.03)},
anchor=south
},
}
\begin{loglogaxis}[
    legend columns = 2,
    legend style = {draw = none},
    scale = 0.75,
    ylabel={$\text{size}_{\text{top}}$},
    xlabel = {$N$},
    ymin=200,
    xtick = \empty,
    xtick = {3*10^4, 1*10^5, 3*10^5, 5*10^5, 8*10^5, 2*10^6},
    xticklabels = {30k, ,0.3M, ,,2M},
    ymajorgrids=true,
    line width=0.25mm,
    grid style=dashed,
    legend entries = { ncols, nrows}
]
\addplot coordinates {
    (32^3, 468) (48^3,658) (64^3,875) (80^3, 1113) (96^3, 1349) (128^3, 1798)
    };
\addplot+ coordinates {
    (32^3, 1274) (48^3,1684) (64^3,2260) (80^3, 2791 ) (96^3, 3333 ) (128^3, 4482)
    };

\addplot [black, domain = 32^3:128^3] {x^(1/3)*10};
\node [ anchor=center] at (0.3*10^6,400) {$\mathcal{O}(N^{1/3})$};
\end{loglogaxis}
\end{tikzpicture}}
    \scalebox{0.7}{
\begin{tikzpicture}
\pgfplotsset{
every axis legend/.append style={
at={(0.5 ,1.03)},
anchor=south
},
}
\begin{loglogaxis}[
    scale = 0.75,
    ylabel={mem$_\text{F}$},
    xlabel = {$N$},
    ymin=200,
    ymin=2*10^6,
    xtick = \empty,
    xtick = {3*10^4, 1*10^5, 3*10^5, 5*10^5, 8*10^5, 2*10^6},
    xticklabels = {30k, ,0.3M, ,,2M},
    ymajorgrids=true,
    line width=0.25mm,
    grid style=dashed,
]
\addplot coordinates {
    (32^3,22067738) (48^3,97924415) (64^3,294666739) (80^3, 623149744) (96^3, 1180405857) (128^3, 3313113881)
    };

\addplot [black, domain = 32^3:128^3] {x*300};
\node [ anchor=center] at (0.3*10^6,3*10^7) {$\mathcal{O}(N)$};
\end{loglogaxis}
\end{tikzpicture}}
    \caption{The growth in the size of the top separator and the memory needed to store the preconditioner for the 3D Inverse Poisson problem with $\alpha \approx 2$. }
    \label{fig:3d_stop_mem}
\end{figure}

In \Cref{fig:3d_stop_mem}, we compare the size of the top separator with the problem size $N$. Note that for a fixed value of $\alpha$, the size of the top separator grows as $\mathcal{O}(N^{1/3})$. Hence, the cost of factorizing the block matrix corresponding to the top separator is $\mathcal{O}(N)$. With $\Theta(\log(N/N_0))$ levels, the total cost of spaQR factorization is $\mathcal{O}(N \log N)$ for a fixed constant aspect ratio $\alpha$. For more detailed complexity analysis, see \Cref{sec:complexity}. The memory needed to store the preconditioner scales as $\mathcal{O}(N)$. Empirically, we see that the scaling is $\mathcal{O}(N^{1.18})$, which may be due to non-asymptotic effects.

\begin{figure}[tbhp]
    \centering
    \scalebox{0.7}{
    \begin{tikzpicture}
    \pgfplotsset{
every axis legend/.append style={
at={(1.02,1)},
anchor=north west},
legend cell align=left}
        \begin{axis}[
        xlabel style = {align=center},
        reverse legend,
            xscale = 1.5,
            yscale = 0.75,
            ybar stacked,
            ymin = 1, ymax = 2000,
            xmin=1, xmax = 15,
            xlabel={Level},
            ylabel={Time (s)},
            xticklabel = \empty,
            x dir=reverse,
            extra x ticks = { 15,14,13,12, 11, 10, 9,8,7,6,5,4,3,2,1 },
            extra x tick labels = {15,,,12,,,9,,,6,,4,,2,},
            enlarge x limits=0.06,
            bar width = 8.5pt,
            ymajorgrids=true,
            line width=0.25mm,
            grid style=dashed,
            xtick align = inside,
            legend image post style={scale=1.7},
            legend entries = { Factorize, Reassign, Scale, Sparsify, Merge}
        ]
        \addplot[ybar, black, pattern=north east lines] table[
                    x = Level,
                    y = Elimination ] {results/runtime_per_level_128_3d_invpoi.dat};
        \addplot+[ybar, black, pattern=grid] table[
                    x = Level,
                    y = Reassign] {results/runtime_per_level_128_3d_invpoi.dat};
        \addplot+[ybar, black, pattern=horizontal lines] table[
                    x = Level,
                    y = Scale] {results/runtime_per_level_128_3d_invpoi.dat};
        \addplot+[ybar, black, pattern=dots] table[
                    x = Level,
                    y = Sparsify] {results/runtime_per_level_128_3d_invpoi.dat};
        \addplot+[ybar, black, pattern=vertical lines] table[
                    x = Level,
                    y = Merge] {results/runtime_per_level_128_3d_invpoi.dat};
        \end{axis}
        \end{tikzpicture}}
    \caption{The runtime per level of the spaQR algorithm split into the five phases: factorize interiors/separators, reassign rows below the diagonal of each factorized interior/separator, scale interfaces, sparsify interfaces, and merge the clusters. We skip sparsification for two levels. The results are shown for the 3D Inverse Poisson problem with $N=128^3$ and $\alpha\approx 2$.}
    \label{fig:runtime_per_level}
\end{figure}
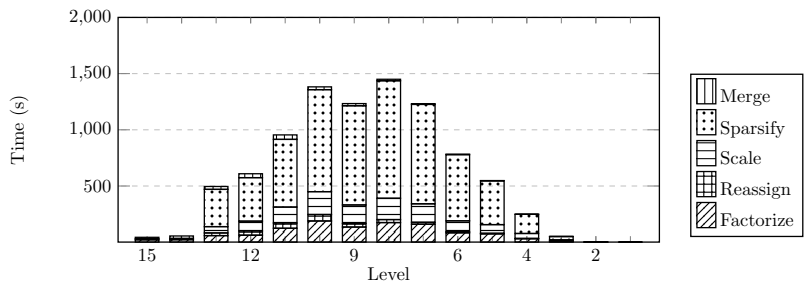

\Cref{fig:runtime_per_level} shows the runtime per level of the spaQR algorithm split into the different phases: factorize interiors/separators, sparsify the extra rows of each interior/separator and reassign them among its ancestors in the elimination tree, scale the interfaces, sparsify the interfaces and merge the clusters. We note that the time per level is roughly the same (a variation of around 6\%) across level 10 to 7. As we scale to bigger problems, the small variations will not matter and we expect a constant runtime per level. 

\subsection{Non-regular problems} Next, we test our algorithm on the set of least-squares problems from the Suite Sparse Matrix Collection~\cite{suitesparse}. The matrix was partitioned using Metis~\cite{metis} and the spaQR algorithm was run with a tolerance $\epsilon$ that is tuned for each problem. Two of the matrices, graphics and Hardesty3, proved the most challenging for spaQR, requiring a lower tolerance of $\epsilon=10^{-4}$ and $\epsilon=10^{-5}$ respectively. CGLS preconditioned with spaQR performed well for all the test matrices and converged in at most 25 iterations as shown in \Cref{Table: spaqr_results}. The time taken to factorize and solve ($t_F + t_S$) with spaQR is much smaller than direct multifrontal QR with the same matrix partitioning. 

\begin{table}[thbp]
    \centering
    
    \caption{Performance of spaQR on some Suite Sparse matrices: tolerance set ($\epsilon$), time to partition, factorize, solve ($t_P$, $t_F$, $t_S$), number of CGLS iterations (\#iter),  memory required to store the preconditioner (mem$_\text{F}$), time to factorize and solve by direct method (Direct)  }
    \scalebox{0.9}{
    \begin{tabular}{@{}lrrrrrrrrc@{}}
        \toprule
        \multicolumn{3}{c}{Matrix} & \multicolumn{6}{c}{spaQR} & Direct \\
       \cmidrule(r){1-3} 
       \cmidrule(lr){4-9} \cmidrule(l){10-10}
        Name & \#rows & \# cols & $\epsilon$ & $t_P$ & $t_F$ & $t_S$ & \#iter  & mem$_\text{F}$ &  $t_F + t_S$\\
        &&& & (s) & (s) & (s) & &  ($10^8$)  & (s)\\
        \midrule
        mesh\_deform & 232k & 9k & $10^{-2}$ & 0.39 & 1.38 & 0.27 & 8  & 0.006 & 17.8 \\
        Kemelmacher & 28k & 9k & $10^{-2}$ & 0.08 & 0.64 & 0.25 & 20 & 
         0.011 & 1.93\\
        human & 20k & 10k & $10^{-2}$ & 0.08  & 0.28  &  0.23 & 25  & 
         0.006 & 1.37 \\
        graphics & 29k & 12k& $10^{-4}$ & 0.07 & 0.25 & 0.13 & 17 & 
         0.005 & 1.19\\
        image\_interp & 240k & 120k & $10^{-3}$ & 0.92 & 5.46 & 1.12 & 8 & 
         0.131 & 73.9 \\
        Hardesty2 & 929k & 303k & $10^{-3}$ & 3.90  & 13.5  & 12.8  & 24  & 
         0.289 & 636 \\
        Hardesty3 & 8.2M & 7.5M & $10^{-5}$ & 96.1  & 384 & 165  & 10  & 
         11.52  & - \\
        \bottomrule
    \end{tabular}}
    \label{Table: spaqr_results}
\end{table}

\section{Conclusions}
\label{sec: conclusions}

In this work, we extended the spaQR algorithm to solve large, sparse linear least squares problems. We introduced a two-step sparsification process to effectively sparsify the interfaces for tall and thin matrices. With this, we empirically showed that the aspect ratio of each diagonal block remains bounded and lead to a $\mathcal{O}(M \log N)$ runtime for the algorithm. We numerically benchmarked our algorithm by testing on large sparse least squares problem arising in optimization problems and other non-regular problems from the Suite Sparse Matrix Collection~\cite{suitesparse}. The numerical tests showed the superiority of spaQR over direct multifrontal QR and CGLS iterative method with a standard diagonal preconditioner.

With this extension to solving large sparse least squares problems, the spaQR solver opens up new applications of low-rank approximations to many areas of computational science. One such area that was investigated in this work is constrained optimization problems. Other potential application domains are computational geometry and computer graphics. This will be investigated in a future work. While the current implementation is sequential, the organization of the algorithm  naturally exhibits parallelism. This will be investigated in future research.

\section{Acknowledgements}
The computing for this project was performed on the Sherlock research cluster, hosted at Stanford University. We thank Stanford University and the Stanford Research Computing Center for providing the computational resources and support that contributed to this research. This work was partly funded by a grant from Sandia National Laboratories (Laboratory Directed Research and Development [LDRD]) entitled ``Hierarchical Low-rank Matrix Factorizations,'' and a grant from the National Aeronautics and Space Administration (NASA, agreement \#80NSSC18M0152). We thank L\'eopold Cambier, Erik G.\ Boman, Juan Alonso, Zan Xu, Jordi Feliu-F\'aba and Steven Brill for valuable discussions. Finally, we also thank Ron Estrin for providing the PDE constrained optimization examples as a test case for our algorithm.

\bibliographystyle{siamplain}
\bibliography{main}

\begin{thebibliography}{10}

\bibitem{mumps}
{\sc P.~R. Amestoy, I.~S. Duff, J.-Y. L'Excellent, and J.~Koster}, {\em Mumps:
  A general purpose distributed memory sparse solver}, in Applied Parallel
  Computing. New Paradigms for HPC in Industry and Academia, T.~S{\o}revik,
  F.~Manne, A.~H. Gebremedhin, and R.~Moe, eds., Berlin, Heidelberg, 2001,
  Springer Berlin Heidelberg, pp.~121--130.

\bibitem{H_QR}
{\sc P.~Benner and T.~Mach}, {\em On the qr decomposition of h-matrices},
  Computing, 88 (2010), \url{https://doi.org/10.1007/s00607-010-0087-y}.

\bibitem{bjorck_article}
{\sc A.~Bjorck}, {\em Stability analysis of the method of semi-normal equations
  for least squares problems}, Linear Algebra and its Applications, 88-89
  (1987), pp.~31--48, \url{https://doi.org/10.1016/0024-3795(87)90101-7}.

\bibitem{Bjorck:1411947}
{\sc A.~Bjorck}, {\em {Numerical methods for least squares problems}}, SIAM,
  Philadelphia, PA, 1996, \url{https://cds.cern.ch/record/1411947}.

\bibitem{2019arXiv190102971C}
{\sc L.~Cambier, C.~Chen, E.~Boman, S.~Rajamanickam, R.~Tuminaro, and
  E.~Darve}, {\em An algebraic sparsified nested dissection algorithm using
  low-rank approximations}, SIAM Journal on Matrix Analysis and Applications,
  41 (2020), pp.~715--746, \url{https://doi.org/10.1137/19M123806X}.

\bibitem{Chow1997ExperimentalSO}
{\sc E.~Chow and Y.~Saad}, {\em Experimental study of ilu preconditioners for
  indefinite matrices}, Journal of Computational and Applied Mathematics, 86
  (1997), pp.~387--414.

\bibitem{suitesparse}
{\sc T.~A. Davis and Y.~Hu}, {\em The university of florida sparse matrix
  collection}, ACM Trans. Math. Softw., 38 (2011),
  \url{https://doi.org/10.1145/2049662.2049663},
  \url{https://doi.org/10.1145/2049662.2049663}.

\bibitem{ZoltanHypergraphIPDPS06}
{\sc K.~D. {Devine}, E.~G. {Boman}, R.~T. {Heaphy}, R.~H. {Bisseling}, and
  U.~V. {Catalyurek}}, {\em Parallel hypergraph partitioning for scientific
  computing}, in Proceedings 20th IEEE International Parallel Distributed
  Processing Symposium, 2006, pp.~10 pp.--,
  \url{https://doi.org/10.1109/IPDPS.2006.1639359}.

\bibitem{estrin2019implementing}
{\sc R.~Estrin, M.~P. Friedlander, D.~Orban, and M.~A. Saunders}, {\em
  Implementing a smooth exact penalty function for equality-constrained
  nonlinear optimization}, SIAM Journal on Scientific Computing, 42 (2020),
  p.~A1809–A1835, \url{https://doi.org/10.1137/19m1238265},
  \url{http://dx.doi.org/10.1137/19M1238265}.

\bibitem{C:LaBRI::CIMI15}
{\sc M.~Faverge, G.~Pichon, P.~Ramet, and J.~Roman}, {\em On the use of
  h-matrix arithmetic in pastix: a preliminary study}, in Workshop on Fast
  Solvers, Toulouse, France, June 2015,
  \url{http://www.labri.fr/~ramet/restricted/cimi15.pdf}.

\bibitem{FeliuFab2018RecursivelyPH}
{\sc J.~Feliu-Fabà, K.~Ho, and L.~Ying}, {\em Recursively preconditioned
  hierarchical interpolative factorization for elliptic partial differential
  equations}, Communications in Mathematical Sciences, 18 (2020), pp.~91--108,
  \url{https://doi.org/10.4310/CMS.2020.v18.n1.a4}.

\bibitem{feliufaba2020hierarchical}
{\sc J.~Feliu-Fabà and L.~Ying}, {\em Hierarchical interpolative factorization
  preconditioner for parabolic equations}, Journal of Scientific Computing, 85
  (2020), \url{https://doi.org/10.1007/s10915-020-01343-5}.

\bibitem{journals/siamsc/FongS11}
{\sc D.~C.-L. Fong and M.~A. Saunders}, {\em Lsmr: An iterative algorithm for
  sparse least-squares problems.}, SIAM J. Scientific Computing, 33 (2011),
  pp.~2950--2971,
  \url{http://dblp.uni-trier.de/db/journals/siamsc/siamsc33.html#FongS11}.

\bibitem{George1973NestedDO}
{\sc A.~{George}}, {\em {Nested Dissection of a Regular Finite Element Mesh}},
  SIAM Journal on Numerical Analysis, 10 (1973), pp.~345--363,
  \url{https://doi.org/10.1137/0710032}.

\bibitem{Ghysels2016AnEM}
{\sc P.~Ghysels, X.~S. Li, F.-H. Rouet, S.~Williams, and A.~Napov}, {\em An
  efficient multicore implementation of a novel hss-structured multifrontal
  solver using randomized sampling}, SIAM Journal on Scientific Computing, 38
  (2016), pp.~S358--S384.

\bibitem{gnanasekaran2020hierarchical}
{\sc A.~Gnanasekaran and E.~Darve}, {\em Hierarchical orthogonal factorization:
  Sparse square matrices}, arXiv preprint arXiv:2010.06807,  (2020).

\bibitem{Golub1965NumericalMF}
{\sc G.~Golub}, {\em Numerical methods for solving linear least squares
  problems}, Numer. Math., 7 (1965), p.~206–216,
  \url{https://doi.org/10.1007/BF01436075},
  \url{https://doi.org/10.1007/BF01436075}.

\bibitem{matrix_computations}
{\sc G.~H. Golub and C.~F. Van~Loan}, {\em Matrix Computations (3rd Ed.)},
  Johns Hopkins University Press, USA, 1996.

\bibitem{FMM_1}
{\sc L.~Greengard and V.~Rokhlin}, {\em A fast algorithm for particle
  simulations}, J. Comput. Phys., 135 (1997), p.~280–292,
  \url{https://doi.org/10.1006/jcph.1997.5706},
  \url{https://doi.org/10.1006/jcph.1997.5706}.

\bibitem{greengard_rokhlin_1997}
{\sc L.~Greengard and V.~Rokhlin}, {\em A new version of the fast multipole
  method for the laplace equation in three dimensions}, Acta Numerica, 6
  (1997), p.~229–269, \url{https://doi.org/10.1017/S0962492900002725}.

\bibitem{Hestenes&Stiefel:1952}
{\sc M.~R. Hestenes and E.~Stiefel}, {\em Methods of conjugate gradients for
  solving linear systems}, Journal of research of the National Bureau of
  Standards, 49 (1952), pp.~409--436.

\bibitem{Ho2016HierarchicalIF}
{\sc K.~Ho and L.~Ying}, {\em Hierarchical interpolative factorization for
  elliptic operators: Integral equations}, Communications on Pure and Applied
  Mathematics, 69 (2013), \url{https://doi.org/10.1002/cpa.21577}.

\bibitem{ho2013hierarchical}
{\sc K.~L. Ho and L.~Ying}, {\em Hierarchical interpolative factorization for
  elliptic operators: Differential equations}, Communications on Pure and
  Applied Mathematics, 69 (2013), pp.~1415--1451.

\bibitem{hsl_mc64}
{\sc {HSL(2013)}}, {\em A collection of fortran codes for large scale
  scientific computation}, \url{http://www.hsl.rl.ac.uk}.

\bibitem{James1990ConjugateGM}
{\sc D.~James}, {\em Conjugate gradient methods for constrained least squares
  problems}, PhD thesis, Dept. of Math., North Carolina State University, 1990.

\bibitem{jennings}
{\sc A.~Jennings and M.~A. Ajiz}, {\em Incomplete methods for solving $a^t ax =
  b$}, SIAM J. Sci. Stat. Comput., 5 (1984), p.~978–987,
  \url{https://doi.org/10.1137/0905067}, \url{https://doi.org/10.1137/0905067}.

\bibitem{metis}
{\sc G.~Karypis and V.~Kumar}, {\em A fast and high quality multilevel scheme
  for partitioning irregular graphs}, SIAM J. Sci. Comput., 20 (1998),
  p.~359–392.

\bibitem{Karypis1998HmetisAH}
{\sc G.~Karypis and V.~Kumar}, {\em Hmetis: a hypergraph partitioning package},
  1998.

\bibitem{klockiewicz2020second}
{\sc B.~Klockiewicz, L.~Cambier, R.~Humble, H.~Tchelepi, and E.~Darve}, {\em
  Second order accurate hierarchical approximate factorization of sparse spd
  matrices}, arXiv preprint arXiv:2007.00789,  (2020).

\bibitem{miqr}
{\sc N.~Li and Y.~Saad}, {\em Miqr: A multilevel incomplete qr preconditioner
  for large sparse least‐squares problems}, SIAM J. Matrix Analysis
  Applications, 28 (2006), pp.~524--550,
  \url{https://doi.org/10.1137/050633032}.

\bibitem{Manteuffel1980AnIF}
{\sc T.~Manteuffel}, {\em An incomplete factorization technique for positive
  definite linear systems}, Mathematics of Computation, 34 (1980),
  pp.~473--497.

\bibitem{FDM_invpoi}
{\sc J.~Nagel}, {\em Numerical solutions to poisson equations using the
  finite-difference method [education column]}, IEEE Antennas and Propagation
  Magazine, 56 (2014), p.~209, \url{https://doi.org/10.1109/MAP.2014.6931698}.

\bibitem{journals/toms/PaigeS82}
{\sc C.~C. Paige and M.~A. Saunders}, {\em Lsqr: An algorithm for sparse linear
  equations and sparse least squares.}, ACM Trans. Math. Softw., 8 (1982),
  pp.~43--71,
  \url{http://dblp.uni-trier.de/db/journals/toms/toms8.html#PaigeS82}.

\bibitem{blr_pastix}
{\sc G.~{Pichon}, E.~{Darve}, M.~{Faverge}, P.~{Ramet}, and J.~{Roman}}, {\em
  Sparse supernodal solver using block low-rank compression}, in 2017 IEEE
  International Parallel and Distributed Processing Symposium Workshops
  (IPDPSW), 2017, pp.~1138--1147.

\bibitem{lorasp1}
{\sc H.~Pouransari, P.~Coulier, and E.~Darve}, {\em Fast hierarchical solvers
  for sparse matrices using extended sparsification and low-rank
  approximation}, SIAM Journal on Scientific Computing, 39 (2017),
  pp.~A797--A830, \url{https://doi.org/10.1137/15M1046939}.

\bibitem{Saad1988PreconditioningTF}
{\sc Y.~Saad}, {\em Preconditioning techniques for nonsymmetric and indefinite
  linear systems}, Journal of computational and applied mathematics, 24 (1988),
  pp.~89--105.

\bibitem{Schmitz2012AFD}
{\sc P.~G. Schmitz and L.~Ying}, {\em A fast direct solver for elliptic
  problems on general meshes in 2d}, J. Comput. Phys., 231 (2012),
  pp.~1314--1338.

\bibitem{iqr_fail}
{\sc X.~Wang}, {\em Incomplete factorization preconditioning for linear least
  squares problems}, PhD thesis, University of Illinois at Urbana-Champaign,
  1994.

\bibitem{Xi2014SuperfastAS}
{\sc Y.~Xi, J.~Xia, S.~Cauley, and V.~Balakrishnan}, {\em Superfast and stable
  structured solvers for toeplitz least squares via randomized sampling}, SIAM
  J. Matrix Analysis Applications, 35 (2014), pp.~44--72.

\bibitem{Xia2013EfficientSM}
{\sc J.~Xia}, {\em Efficient structured multifrontal factorization for general
  large sparse matrices}, SIAM J. Scientific Computing, 35 (2013).

\bibitem{Xia2009SuperfastMM}
{\sc J.~Xia, S.~Chandrasekaran, M.~Gu, and X.~S. Li}, {\em Superfast
  multifrontal method for large structured linear systems of equations}, SIAM
  J. Matrix Analysis Applications, 31 (2009), pp.~1382--1411.

\bibitem{lorasp2}
{\sc K.~Yang, H.~Pouransari, and E.~Darve}, {\em Sparse hierarchical solvers
  with guaranteed convergence}, International Journal for Numerical Methods in
  Engineering,  (2016), \url{https://doi.org/10.1002/nme.6166}.

\bibitem{atalyrek2011PaToHT}
{\sc {\"U}.~V. Çataly{\"u}rek and C.~Aykanat}, {\em Patoh (partitioning tool
  for hypergraphs)}, in Encyclopedia of Parallel Computing, 2011.

\end{thebibliography}

\appendix
\section{Matrix Generation: 2D Inverse Poisson problem}
\label{sec: mat_gen_2d}
    
\begin{figure}
    \centering
    \scalebox{0.65}{
      \begin{tikzpicture}
    \foreach \x in {0,2,...,4}{
        \foreach \y in {0,2,...,4}{
            \fill (\x,\y) circle (3pt);
        }
    }
    \draw[step=2cm] (0,0) grid (4,4);
    \foreach \x in {1,3}{
        \foreach \y in {1,3}{
            \draw (\x,\y) node[cross] {};
        }
    }
    \node at (-0.5,-0.4) {\large$u_{i-1,j-1}$};
    \node at (2,-0.4) {\large$u_{i,j-1}$};
    \node at (4.5,-0.4) {\large$u_{i+1,j-1}$};
    \node at (-0.5,2.25) {\large$u_{i-1,j}$};
    \node at (2.5,2.25) {\large$u_{i,j}$};
    \node at (4.6,2.25) {\large$u_{i+1,j}$};
    \node at (-0.5,4.4) {\large$u_{i-1,j+1}$};
    \node at (2,4.4) {\large$u_{i,j+1}$};
    \node at (4.5,4.4) {\large$u_{i+1,j+1}$};
    \node at (1,1.4) {\large$z_{i-1,j-1}$};
    \node at (3,1.4) {\large$z_{i,j-1}$};
    
    \node at (1,3.4) {\large$z_{i-1,j}$};
    \node at (3,3.4) {\large$z_{i,j}$};
  \end{tikzpicture}}
  \caption{Finite difference mesh for the 2D variable coefficient Poisson equation.}
  \label{fig: stag_grid}
\end{figure}
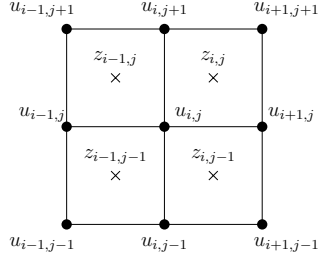

    Consider the variable coefficient Poisson equation,
    \begin{align*}
      -\nabla \cdot (z \nabla u) &= h \text{ in }\Omega  \\
      u = 0 \text{ in }d\Omega
    \end{align*}
    discretized using finite difference on a staggered grid~\cite{FDM_invpoi} shown in \Cref{fig: stag_grid}. The discretization at grid point $(i,j)$ is, 
    \begin{equation} \label{eq:poisson}
    f(u, z) \coloneqq -a_0u_{i,j}+a_1u_{i+1,j}+a_2u_{i,j+1}+a_3u_{i-1,j}+a_4u_{i,j-1} + q_{i,j} = 0
\end{equation}
    \begin{align*}
    a_0 &= z_{i,j}+z_{i-1,j}+z_{i,j-1}+z_{i-1,j-1} & \\
    a_1 &= \frac{1}{2}(z_{i,j}+z_{i,j-1}) \qquad  &a_2 &= \frac{1}{2}(z_{i-1,j}+z_{i,j})\\
    a_3 &= \frac{1}{2}(z_{i-1,j-1}+z_{i-1,j}) \qquad &a_4 &= \frac{1}{2}(z_{i,j-1}+z_{i-1,j-1}) \\
    q_{i,j} &= \int_{\Omega_{ij}} \!\! h(r) \, d\Omega \text{,} \quad \Omega_{ij} \text{ is a square region around } u_{ij} 
\end{align*}

The finite difference discretization leads to $n^2$ equations, one for each grid point on the 2D $n \times n$ grid. The number of variables is $n_u +n_z$, where $n_u = n^2$ and $n_z = (n+1)^2$. The Jacobian matrix is constructed by taking partial derivatives of each of the $n^2$ equations of the form \Cref{eq:poisson}, with respect to each $u(i,j)$ and $z(i,j)$ variable. For example, the following are the partial derivatives of \Cref{eq:poisson},
\[
    \frac{\partial f}{\partial u_{i,j}} = -a_0  \quad
    \frac{\partial f}{\partial u_{i+1,j}} = a_1  \quad
    \frac{\partial f}{\partial u_{i,j+1}} = a_2  \quad
    \frac{\partial f}{\partial u_{i-1,j}} = a_3  \quad
    \frac{\partial f}{\partial u_{i,j-1}} = a_4
\]
\begin{align*}
    \frac{\partial f}{\partial z_{i,j}} &= -u_{i,j}+\frac{1}{2}u_{i+1,j}+ \frac{1}{2}u_{i,j+1} \qquad &\frac{\partial f}{\partial z_{i-1,j}} &= -u_{i,j}+ \frac{1}{2}u_{i,j+1}+ \frac{1}{2}u_{i-1,j} \\
    \frac{\partial f}{\partial z_{i,j-1}} &= -u_{i,j}+\frac{1}{2}u_{i+1,j}+ \frac{1}{2}u_{i,j-1} \qquad 
    &\frac{\partial f}{\partial z_{i-1,j-1}} &= -u_{i,j}+ \frac{1}{2}u_{i-1,j}+ \frac{1}{2}u_{i,j-1} 
\end{align*}

We can evaluate the Jacobian at random values of $u, z$. Then, we solve the least squares problem, $J^T x = b$. In general, choosing random values of $u, z$ leads to $J^T$ having an aspect ratio of 2. It is possible to reduce the aspect ratio of the matrix by artificially creating rows that are equal to 0. This is done by choosing values for $u$ such that $u$ is constant in some region, say $u=1$ in some region of the mesh. In that case, \Cref{eq:poisson} simplifies and the entire row become zero. Then, by removing zero rows, we can reduce the aspect ratio of matrix $J^T$. 

For example, setting $u=1$ and $z=1$
everywhere on the grid (except the boundary where $u=0$) will make most of the partial derivatives with respect to $z$ variables zero (expect the $z$ variables near the boundary of the mesh). This leads to zero rows in $J^T$ which can be removed to get an aspect ratio close to 1 for $J^T$. 

\section{Matrix Generation: 3D Inverse Poisson problem}
\label{sec: mat_gen_3d}

The finite difference discretization for the variable coefficient Poisson equation in a 3D $n \times n \times n$ mesh is as follows. At grid point $(i,j)$, 
\begin{align}
    f(u,z) &\coloneqq -a_0 u_{i,j,k} + a_1 u_{i+1,j,k} + a_2 u_{i,j+1,k} + a_3 u_{i,j,k+1} \\
     & \quad + a_4 u_{i-1, j, k} + a_5 u_{i, j-1, k} + a_6 u_{i, j, k-1} + q_{i,j,k} = 0 \nonumber
\end{align}
\begin{align*}
    a_0 &= \frac{3}{4}\big( z_{i,j,k} + z_{i-1,j,k} + z_{i, j-1, k} + z_{i, j, k-1} \\ 
    & \quad \quad  + z_{i-1, j-1, k} + z_{i, j-1, k-1} + z_{i-1, j, k-1} + z_{i-1, j-1, k-1}  \big) \\
    a_1 &= \frac{-1}{4}\big( z_{i,j,k} + z_{i,j, k-1} + z_{i, j-1, k} + z_{i, j-1, k-1} \big) \\
    a_2 &= \frac{-1}{4}\big( z_{i-1,j,k} + z_{i,j, k} + z_{i-1, j, k-1} + z_{i, j, k-1} \big) \\
    a_3 &= \frac{-1}{4}\big( z_{i-1,j,k} + z_{i,j, k} + z_{i-1, j-1, k} + z_{i, j-1, k} \big) \\
    a_4 &= \frac{-1}{4}\big( z_{i-1,j,k} + z_{i-1,j, k-1} + z_{i-1, j-1, k} + z_{i-1, j-1, k-1} \big) \\
    a_5 &= \frac{-1}{4}\big( z_{i-1,j-1,k} + z_{i,j-1, k} + z_{i-1, j-1, k-1} + z_{i, j-1, k-1} \big) \\
    a_6 &= \frac{-1}{4}\big( z_{i-1,j,k-1} + z_{i,j, k-1} + z_{i-1, j-1, k-1} + z_{i, j-1, k-1} \big) \\
    q_{i,j,k} &= \int_{\Omega_{ijk}} \!\! h(r) \, d\Omega \text{,} \quad \Omega_{ijk} \text{ is a cubic region around } u_{i,j,k}
\end{align*}

For the 3D problem, we have $n^3$ equations, one for each grid point on the $n \times n \times n$ mesh. There are $n_u + n_z$ variables, where $n_u = n^3$ and $n_z = (n+1)^3$. The Jacobian can be constructed by taking partial derivatives of the discretized equations with respect to each $u$ and $z$ variable. As with the 2D problem, by evaluating the Jacobian for specific values of $u$ and $z$, we can generate least squares matrices ($J^T$) with varying aspect ratio. 

\end{document}